\documentclass[12pt,draftcls,onecolumn]{IEEEtran}
\usepackage{amssymb}
\usepackage[cmex10]{amsmath}
\usepackage{amsfonts}
\usepackage[latin1]{inputenc}
\usepackage{makeidx}         
\usepackage{placeins}
\usepackage{pdfsync}
\usepackage{multicol}        
\usepackage[bottom]{footmisc}
\usepackage{graphicx}
\usepackage{varioref}
\usepackage{url}
\usepackage{cancel}

\usepackage{placeins}
\usepackage{color}

\newcommand{\tp}{^\top}

\newcommand{\E}[1]           {\mathop{\rm\,E} \left[ #1 \right]}        
\newcommand{\LL}{{\mathcal{L}_+^{\Gamma}}}                                                           
\newcommand{\BL}{{\partial{\mathcal{L}_+^{\Gamma}}}}                                      
\newcommand{\CL}{{\overline{\mathcal{L}_+^{\Gamma}}}}                                 

\newcommand{\tr}{\mathop{\rm tr}}                                                                            

\newcommand{\imunit}{{\rm j}}
\newcommand{\e}{{\rm\,e}}

\newcommand{\Range}{\mathop{\rm Range}}

\newcommand{\eq}{\begin{equation}}
\newcommand{\eeq}{\end{equation}}
\newcommand{\eqn}{\begin{eqnarray}}
\newcommand{\eeqn}{\end{eqnarray}}

\newcommand{\bsea}{\begin{subeqnarray}}
\newcommand{\esea}{\end{subeqnarray}}

\newcommand{\symm}{\mathcal{Q}}
\newcommand{\DDdue}[2]{\mathbb{D}\left(#1\|#2\right)}
\newcommand{\be}{\begin{equation}}
\newcommand{\ee}{\end{equation}}
\newcommand{\bbm}{\begin{bmatrix}}
\newcommand{\ebm}{\end{bmatrix}}

\newcommand{\Ct}{\mathcal{V}\left(\mathcal{C}_{R+}^{m\times m}\right)}
\newcommand{\Ctpiu}{\mathcal{C}_{R+}^{m\times m}}

\newfont{\BB}{msbm10 scaled\magstep1}
\newfont{\bb}{msbm8}

\def\Z{\mbox{\BB Z}}

\def\D{\mbox{\BB D}}
\def\E{\mbox{\BB E}}

\def\T{\mathbb{T}}

\def\epi{{\rm epi}\ }

\newcommand{\Rgamma}{\mathop{\rm Range} \Gamma}

\newcommand{\Jpsi}{{J_\Psi}}

\newcommand{\R}{\mathbb{R}}

\newcommand{\C}{\mathbb{C}}

\newcommand{\Specspacepiu}{\mathcal{S}_+^{m\times m}}

\newcommand{\ejth}{{\rm e}^{ {\imunit}\vartheta}}
\newcommand{\suchthat}{\mbox{such that}}

\newcommand{\trasp}{^\top}

\newcommand{\dF}{dF(\vartheta)}
\newcommand{\h}[1]{h_r({#1})}
\newcommand{\dr}[2]{\D_r(#1\|#2)}
\newcommand{\deltayk}{\delta \hat{y}_k}
\newcommand{\deltayka}[1]{\delta \hat{y}_{#1}}
\newcommand{\deltazk}{\delta \hat{z}_k}
\newcommand{\deltazka}[1]{\delta \hat{z}_{#1}}

 \definecolor{darkGreen}{rgb}{0,0.6,0}

\definecolor{Royalblue}{cmyk}{1,0.30,0.2,0.2}

\DeclareMathOperator*{\argmin}{arg\,min}
\def\bmat{\left[ \begin{array}}
\def\emat{\end{array} \right]}

\newcounter{pippo}

\newtheorem{remark}{Remark}[section]
\newtheorem{theorem}{Theorem}[section]
\newtheorem{corr}{Corollary}[section]
\newtheorem{propo}{Proposition}[section]
\newtheorem{lemm}{Lemma}[section]
\newtheorem{exam}{Example}
\newtheorem{probl}[pippo]{Problem}
\newtheorem{definition}{Definition}[section]
\newtheorem{apropo}{Proposition}[subsection]
\newtheorem{alemm}{Lemma}[subsection]

\newcommand{\teo}{\begin{teor}}
\newcommand{\eteo}{\end{teor}}
\newcommand{\cor}{\begin{corr}}
\newcommand{\ecor}{\end{corr}}
\newcommand{\prop}{\begin{propo}}
\newcommand{\eprop}{\end{propo}}
\newcommand{\lem}{\begin{lemm}}
\newcommand{\elem}{\end{lemm}}
\newcommand{\ex}{\begin{exam}}
\newcommand{\eex}{\end{exam}}
\newcommand{\pb}{\begin{probl}}
\newcommand{\epb}{\end{probl}}
\newcommand{\df}{\begin{defn}}
\newcommand{\edf}{\end{defn}}
\newcommand{\aprop}{\begin{apropo}}
\newcommand{\eaprop}{\end{apropo}}
\newcommand{\alem}{\begin{alemm}}
\newcommand{\ealem}{\end{alemm}}

\makeindex
\begin{document}



\title{
Time and spectral domain relative entropy: \\A new approach to multivariate spectral estimation}


\author{Augusto~Ferrante, Chiara~Masiero and~Michele~Pavon\thanks{Work partially supported by the Italian Ministry for Education and Research (MIUR) under PRIN grant n. 20085FFJ2Z
``New Algorithms and Applications of System Identification and Adaptive Control"}
\thanks{A. Ferrante and C. Masiero are  with the
Dipartimento di Ingegneria dell'Informazione, Universit\`a di
Padova, via Gradenigo 6/B, 35131 Padova, Italy {\tt\small
augusto@dei.unipd.it},  {\tt\small
masiero.chiara@dei.unipd.it}}\thanks{M. Pavon is with the Dipartimento
di Matematica Pura ed Applicata, Universit\`a di Padova, via
Trieste 63, 35131 Padova, Italy {\tt\small pavon@math.unipd.it}}}

\markboth{DRAFT}{Shell \MakeLowercase{\textit{et al.}}: Bare Demo of IEEEtran.cls for Journals}

\maketitle

\begin{abstract}
The concept of {\em spectral relative entropy rate} is introduced for jointly stationary Gaussian processes. Using
{classical information-theoretic results}, we establish a remarkable connection between time and spectral domain relative entropy rates.
This naturally leads to a new  { spectral estimation technique where a multivariate version of the {\em Itakura-Saito distance} is employed}. { It} may be viewed as an extension of the approach, called THREE, introduced by Byrnes, Georgiou and Lindquist in 2000 which, in turn, followed in the footsteps of the Burg-Jaynes Maximum Entropy Method.  Spectral estimation is here recast in the form of a constrained spectrum approximation problem where the distance is equal to the processes relative entropy rate. The corresponding solution entails a complexity upper bound which improves on the one so far available in the multichannel framework. Indeed, it is {\em equal} to the one featured by THREE in the scalar case. The solution is computed via a globally convergent matricial Newton-type algorithm. Simulations suggest the effectiveness of the new technique in tackling multivariate spectral estimation tasks, especially in the case of short data records.

\end{abstract}


\title{
 A novel maximum entropy problem to ensure feasibility of
THREE-like high-resolution spectral estimators}


\section{Introduction} \label{cap1}
Multidimensional spectral estimation is an old and challenging problem \cite{McClellan_Multidimensionalspectralestimation,STOICA_MOSES_SPECTRALANALYSIS} which keeps generating widespread interest in the natural and engineering sciences, see e.g. \cite{GEORGIOU_MAXIMUMENTROPY,GEORGIOU_RELATIVEENTROPY,ROSEN_STOFFER_MULTIVARIATESPECTRA,RAMPONI_FERRANTE_PAVON_GLOBALLYCONVERGENT}. A new approach to scalar spectral estimation  called THREE was introduced by Byrnes, Georgiou and
Lindquist in  \cite{BYRNES_GEORGIOU_LINDQUIST_THREE,GEORGIOU_SELECTIVEHARMONIC}. It may be viewed as a (considerable) generalization of  classical Burg-like maximum entropy methods. This estimator
permits higher resolution in prescribed frequency bands and
is particularly competitive in the case of short
observation records. In this approach, the output covariance of a bank of filters  is used to extract information on the input power spectrum. A first attempt to generalize this approach to the multichannel situation was made in \cite{RAMPONI_FERRANTE_PAVON_GLOBALLYCONVERGENT}, where, due to the lack of multidimensional theoretical results, a non entropy-like distance was employed in the optimization part of the procedure. The resulting solution, however, had higher McMillan degree than in the original scalar THREE method.

The main contribution of this paper is twofold: On the one hand, we introduce what appears to be a  most natural multivariate { generalization} of the THREE method, called RER { since} the metric employed in the { optimization problem} originates from the {\em relative entropy rate} of the two processes.  { The latter may be viewed} {  as a multivariate extension of the classical {\em Itakura-Saito distance} { widely used in} signal processing \cite{BASSEVILLE_DISTANCEMEASURES}}. { The proposed} method features a complexity upper bound which, considerably improving on the one so far available, is in fact {\em equal} to the one featured by THREE in the scalar case. Like all previous THREE-like methods, RER exhibits high resolution features and works extremely well, outperforming classical identification methods, in the case of short observation records.
 On the other hand, further cogent support for the choice of our distance measure between spectra is provided by  a novel information-theoretic { result: We  introduce the concept of {\em spectral entropy rate}  for stationary Gaussian processes and we establish a {\em circular symmetry property} of  the increments of the process occurring in the spectral representation. Then,} using classical results of  Pinsker \cite{PINSKER_INFORMATION}, Van den Bos \cite{van_den_bos95}, Stoorvogel and Van Schuppen \cite{STOORVOGEL_VANSCHUPPEN_SYSTEMCRITERIA}, we prove that the time and spectral domain relative entropy rates are in fact equal! This profound result is deferred to the last section of the paper for expository reasons.

The paper is outlined as follows.
Section \ref{BACK} collects basic results on entropy for Gaussian vectors and processes.  { Section \ref{introTHREE} introduces THREE-like spectral estimation methods. Section \ref{nnsec4} presents the new approach RER via a convex optimization problem and derives the form of the optimal spectral estimate. In Section \ref{dualproblem}, we establish a nontrivial existence result for the dual problem. A globally convergent, matricial  Newton-type method is presented in Section \ref{nnsec6} to solve the dual problem. The computational burden is dramatically reduced
thanks to various nontrivial results of {\em spectral factorization}.} In Section \ref{simul}, both scalar and multivariate examples are studied via simulation: the performance of the RER method { is} compared to { that of previously available approaches}. In Section \ref{spectralprelim}, some background results on complex Gaussian random vectors and on the spectral representation of stationary Gaussian processes are presented.
Finally, in Section \ref{SRER},  we introduce the spectral relative entropy rate of Gaussian processes and { establish a profound} connection between time and spectral domain relative entropy rates.

\section{Background on entropy for Gaussian processes}\label{BACK}

{ We collect below some basic concepts} and results on entropy of  Gaussian random vectors and processes that may be found e.g. in \cite{PINSKER_INFORMATION,IHARA_INFORMATION,COVER_THOMAS}. The {\em differential entropy} $H(p)$ of a probability density function $p$ on $\R^n$ is defined by
\begin{equation}\label{DiffEntropy}
H(p)=-\int_{\R^n}\log (p(x))p(x)dx.
\end{equation}
In case of a zero-mean Gaussian { density} $p$ with nonsingular covariance matrix $P$, we get
\begin{equation}\label{gaussianentropy}
 H(p)=\frac{1}{2}\log(\det P)+\frac{1}{2}n\left(1+\log(2\pi)\right),
\end{equation}
The {\em relative entropy} or { {\em Kullback-Leibler}} pseudo-distance or {\em divergence} between two probability densities $p$ and $q$, with the support of $p$ contained in the support of $q$, is defined by
  \begin{equation}\D(p\|q){ :=}\int_{\R^n} p(x)\log\frac{p(x)}{q(x)}dx,
\end{equation}
see e.g \cite{COVER_THOMAS}.
In the case of two zero-mean Gaussian { densities} $p$  and $q$ with { positive definite} covariance  matrices  $P$ and $Q$, respectively, the relative entropy is given by:
\begin{equation}\label{divgauss}\D(p\|q)=\frac{1}{2}\left[ \log\det
(P^{-1}Q)+\tr(Q^{-1}P)-n \right].
\end{equation}
Consider now a discrete-time { Gaussian} process $\{y_k;\,k\in\Z\}$ taking values in $\R^m$. Let  $Y_{[-n,n]}$ be the random vector obtained by considering the window $y_{-n},y_{-n+1},\cdots,y_0,\cdots,$ $y_{n-1},y_n$, and let $p_{Y_{[-n,n]}}$ denote the corresponding joint density.
\begin{definition} The {\em (differential) entropy rate} of $y$ is defined by
\begin{equation}\label{ER}
{\h{y} :=}\lim_{n\rightarrow\infty}\frac{1}{2n+1}H(p_{Y_{[-n,n]}}),
\end{equation}
if the limit exists.
\end{definition}
In \cite{KOL56}, Kolmogorov established the following important result.
\begin{theorem} Let $y=\{y_k;\,k\in\Z\}$ be a $\R^m$-valued, zero-mean, Gaussian, stationary, purely nondeterministic stochastic process with spectral density $\Phi_y$.
Then
\begin{equation}\label{entropyconnection}
\h{y}=\frac{m}{2}\log (2\pi e)+\frac{1}{4\pi}\int_{-\pi}^\pi\log\det\Phi_y(e^{{\imunit}\vartheta})d\vartheta.
\end{equation}
\end{theorem}
As is well-known, there is also a fundamental connection between the quantity appearing in
(\ref{entropyconnection}) and the optimal one-step-ahead predictor: The multivariate Szeg\"{o}-Kolmogorov formula. It is
\begin{equation}\label{SK}
\det R=\exp\left\{\frac{1}{2\pi}\int_{-\pi}^\pi\log\det\Phi_y(e^{{\imunit}\vartheta})d\vartheta\right\},
\end{equation}
where $R$ is the error covariance matrix corresponding to the optimal predictor.
Let $y=\{y_k;\,k\in\Z\}$, $z=\{z_k;\,k\in\Z\}$ be two zero-mean, jointly Gaussian, stationary, purely nondeterministic  processes taking values in $\R^m$.
Let $Y_{[-n,n]}$ and $Z_{[-n,n]}$ be defined as above.
\begin{definition} The {\em relative entropy rate} between $y$ and $z$ is defined by
\begin{equation}\label{DR}
{\dr{y}{z} :=} \lim_{n\rightarrow\infty}\frac{1}{2n+1}\D(p_{Y_{[-n,n]}}\|p_{Z_{[-n,n]}})
\end{equation}
if the limit exists.
\end{definition}
The following interesting result {holds (see \cite{STOORVOGEL_VANSCHUPPEN_SYSTEMCRITERIA,IHARA_INFORMATION})}.
\begin{theorem}Let $y=\{y_k;\,k\in\mathbb{Z}\}$ and $z=\{z_k;\,k\in\mathbb{Z}\}$ be $\R^m$-valued, zero-mean, Gaussian, stationary, purely  nondeterministic processes with spectral density functions $\Phi_y$ and $\Phi_z$, respectively.
Assume, moreover, that  at least one of the following conditions is satisfied:
\begin{enumerate}
\item $\Phi_y\Phi_z^{-1}$ is bounded;
\item $\Phi_y\in\mathrm{L}^{2}\left(-\pi,\pi\right)$ and $\Phi_z$ is {\em coercive} (i.e. $\exists\,\alpha>0$ s.t.  $\Phi_z(e^{{\imunit}\vartheta})-\alpha I_m>0$ a.e. on $\T$).
\end{enumerate}
Then
\begin{equation}\label{relentropyconnection}
\dr{y}{z}=\frac{1}{4\pi}\int_{-\pi}^\pi\left\{\log\det \left(\Phi_y^{-1}(e^{{\imunit}\vartheta})\Phi_z(e^{{\imunit}\vartheta})\right)+\tr\left[\Phi_z^{-1}(e^{{\imunit}\vartheta})\left(\Phi_y(e^{{\imunit}\vartheta})-\Phi_z(e^{{\imunit}\vartheta})\right)\right]\right\}d\vartheta.
\end{equation}
\end{theorem}

\section{THREE-like estimation and generalized moment problems}\label{introTHREE}
We denote by
 $\mathcal{S}_+ ^{m\times m}$  the family of bounded and coercive spectral densities on $\T:=\{z\in\C:\ |z|=1\}$ of $\R^m$-valued processes.
Suppose that the data $\{y_i\}_{i=1}^N$ are generated by an unknown, zero-mean, $m$-dimensional,  $\R^m$-valued, purely nondeterministic, stationary, Gaussian process $y=\{y_k;\, k\in\,\mathbb{Z}\}$.  We wish to estimate the spectral density $\Phi \in \, \mathcal{S}_+ ^{m\times m}$ of $y$ from $\{y_i\}_{i=1}^N$.  A THREE-like approach generalizes Burg-like methods in several ways. The second order statistics that are estimated from the data $\{y_i\}_{i=1}^N$ are not necessarily the covariance lags $C_l:=E\{y_{k+l} y_k\trasp\}$ of $y$. Moreover, a prior estimate of $\Phi$ may be included in the estimation procedure. More explicitly these methods hinge on the following four elements:
\begin{enumerate}
 \item A rational filter to process the data. The filter has transfer function
\begin{equation}\label{eq:filter}
 G(z)={(zI-A)}^{-1}B,
\end{equation}
where $A\in\mathbb{R}^{n\times n}$ is a stability matrix (i.e. it has all its eigenvalues inside the unit circle), $B\in \mathbb{R}^{n\times m}$ is full rank, $n\geq m$, and $(A,B)$ is a reachable pair;
 \item an estimate based on the data $\{y_i\}_{i=1}^N$ of the steady-state covariance $\Sigma$ of the state $x(k)$ of the filter
 \begin{equation}\label{statefilter}
 x(k+1)=Ax(k)+By(k);
 \end{equation}
 \item a \emph{prior} spectral density $\Psi\in \mathcal{S}_+ ^{m\times m}$;
 \item an index that measures the distance between two spectral densities.
\end{enumerate}
The filterbank (\ref{statefilter}) provides Carath\`eodory or, more generally, Nevanlinna-Pick interpolation data for  the positive real part $\Phi_+$ of $\Phi$,
see \cite[Section II]{BYRNES_GEORGIOU_LINDQUIST_THREE}. This occurs  through the constraint
\begin{equation}\label{CONST}
\int G\Phi G^*=\Sigma
\end{equation}
which must be satisfied by the spectrum of $y$ (here and throughout the paper, integration | when not otherwise specified | is on the unit circle with respect to normalized Lebesgue measure). Concerning the spectral density $\Psi$:  It allows to take into account possible \emph{a priori} information on $\Phi$, a contingency that is frequent in practice. For example, $\Psi$ may simply be a coarse estimate of the true spectrum.\footnote{When no prior information on $\Phi$ is available, $\Psi$ is set   either to the identity or to the sample covariance of the available data $\{y_i\}_{i=1}^N$. Dually, the prior $\Psi$ yields a smooth parameterization of solutions with bounded degree which permits tuning.} Since, in general, $\Psi$ is not consistent with the interpolation conditions, an approximation problem arises. It is then necessary to introduce an adequate distance index. This crucial choice is dictated by several requirements. On the one hand, the solution should be rationa!
 l of low McMillan degree at least when the prior $\Psi$ is such. On the other hand, the variational analysis should lead to a computable solution, typically by solving the dual optimization problem. In the scalar case \cite{BYRNES_GEORGIOU_LINDQUIST_THREE,GEORGIOU_LINDQUIST_KULLBACKLEIBLER}, the choice { was made of minimizing the following  Kullback-Leibler type criterion:
$d_{KL}(\Psi,\Phi) = \int \Psi \log \frac{\Psi}{\Phi}.$}
This choice features both of the above specifications. In the multivariable case, a Kullback-Leibler pseudo-distance
may also be readily defined \cite {GEORGIOU_RELATIVEENTROPY}, inspired by the {\em Umegaki-von Neumann's relative entropy}  \cite{NIELSEN_CHUANG_QUANTUMCOMPUTATION} of statistical quantum mechanics. The resulting  spectrum approximation problem, however, leads to computable solutions of bounded McMillan degree only in the case when the prior spectral density has the form $\Psi(z)=\psi(z) I$, where $\psi(z)$ is a scalar spectral density (yielding the {\em maximum entropy solution} when $\Psi=I$, \cite{GEORGIOU_MAXIMUMENTROPY,BLOMQVIST_LINDQUIST_NAGAMUNE_MATRIXVALUED,GEORGIOU_RELATIVEENTROPY}).
On the contrary, with the following multivariate extension of the  Hellinger distance introduced in \cite{FERRANTE_PAVON_RAMPONI_HELLINGERVSKULLBACK},
\begin{equation}
\label{HELLINGER}
\begin{split}
d_H(\Psi,\Phi)^2 &:= \inf_{W_\Psi, W_\Phi} \tr \int \left(W_\Psi-W_\Phi\right)\left(W_\Psi-W_\Phi\right)^\ast, \\
&\suchthat \quad W_\Psi W_\Psi^\ast = \Psi \quad {\rm and} \quad  W_\Phi W_\Phi^\ast = \Phi,
\end{split}
\end{equation}
which is a {\em bona fide} distance, the variational analysis can be carried out leading to a computable solution ((\ref{HELLINGER}) is just the $L^2$-distance between the sets of {\em square spectral factors} of the two spectra). An effective multivariate THREE-like spectral estimation method can then be based on such a distance, leading to rational solutions when the prior $\Psi$ is rational \cite{RAMPONI_FERRANTE_PAVON_GLOBALLYCONVERGENT}. The complexity of the solution, however, is usually { noticeably} higher than in the original scalar THREE approach.\\
We show that, employing the relative entropy rate  \eqref{relentropyconnection} as  index for the approximation problem, the variational analysis can be carried out explicitly. { Moreover, }such a choice yields an upper bound on the complexity of the solution equal to that in the original THREE method.
\begin{remark}
Notice that finding an input process that is compatible with the estimated covariance and has {\em rational} spectrum of prescribed maximum degree turns into a Nevanlinna-Pick interpolation problem with bounded degree \cite{BLOMQVIST_LINDQUIST_NAGAMUNE_MATRIXVALUED},\cite{GEORGIOU_INTERPOLATION}.
The latter can be viewed as a {\em generalized moment problem} which is advantageously cast in the frame of various convex optimization problems. An example is provided by the {\em covariance extension problem} and its generalization, see
\cite{GEORGIOU_REALIZATION},
\cite{COMPLETE_PARAMETERIZATION}
\cite{PARTIAL_STOCHASTIC_REALIZATION}
\cite{BYRNES_GUSEV_LINDQUIST_CONVEX},
\cite{BYRNES_GEORGIOU_LINDQUIST_GENERALIZEDCRITERION},
\cite{GEORGIOU_MAXIMUMENTROPY}.
These problems pose a number of theoretical and computational challenges for which we also refer the reader to
\cite{GEORGIOU_LINDQUIST_KULLBACKLEIBLER},
\cite{GEORGIOU_STATECOV},
\cite{GEORGIOU_GENERALMOMENT}, and
\cite{BYRNES_LINDQUIST_IMPORTANTMOMENTS}.
Besides signal processing, significant applications of this theory are found in modeling and identification
\cite{BYRNES_ENQVIST_LINDQUIST_IDENTIFIABILITY},
\cite{GEORGIOU_LINDQUIST_CONVEXARMA},
\cite{Enqvist_Karlsson_covariance-Interpolation},
$H_\infty$ robust control
\cite{BYRNES_GEORGIOU_LINDQUIST_MEGRETSKI},
\cite{GEORGIOU_LINDQUIST_REMARKSDESIGN},
and biomedical engineering
\cite{NASIRI_EBBINI_GEORGIOU_NONINVASIVEESTIMATION}.
\end{remark}
\begin{remark}
 In spectral estimation, it is important to develop problem-specific criteria for choosing a spectral density from a given family  satisfying prescribed constraints or to be able to compare such spectral densities in an informative, quantitative manner. For instance, in \cite{GEORGIOU_DISTANCES, GEORGIOU_RIEMANNIAN,GEORGIOU_JIANG_NING}, it was shown that a geometry entirely analogous to the geometry of the Fisher information metric exists for power spectral densities. Moreover, distances between power spectra can be used quite effectively in identifying transitions, changes, and affinity between time series or even spacial series.
Applications include automated phoneme recognition by identifying natural transition time markers in speech, which separate segments of maximal spectral separation using a suitable metric and variants thereof. Along a similar line, two-dimensional distributions are identified on the inside and outside of a curve, and then, the curve is evolved using geometric active contours to ensure maximal separation of the spectral content of two regions. This idea has been recently applied to visual tracking  \cite{GEORGIOU_MEANING_DISTANCE}.

 \end{remark}

\section{A new metric for multivariate spectral estimation}
\label{nnsec4}

Motivated by relation (\ref{relentropyconnection}), we define a new pseudo-distance among spectra in $\mathcal{S}_+ ^{m\times m}$:
\begin{equation}\label{RERDIST}
d_{RER}(\Phi,\Psi):=\frac{1}{4\pi}\int_{-\pi}^\pi\left\{\log\det \left(\Phi^{-1}(e^{{\imunit}\vartheta})\Psi(e^{{\imunit}\vartheta})\right)+\tr\left[\Psi^{-1}(e^{{\imunit}\vartheta})\left(\Phi(e^{{\imunit}\vartheta})-\Psi(e^{{\imunit}\vartheta})\right)\right]\right\}d\vartheta.
\end{equation}
{ Further motivation for this distance choice is provided by a profound, information-theoretic result relating time and spectral domain relative entropy rates, see Theorem \ref{equalrer} below. Notice that in the case of scalar spectra, $d_{RER}(\Phi,\Psi)=1/2 d_{IS}(\Phi,\Psi)$, where
$$d_{IS}(\Phi,\Psi)=\frac{1}{2\pi}\int_{-\pi}^\pi\left\{\frac{\Phi(e^{{\imunit}\vartheta})}{\Psi(e^{{\imunit}\vartheta})}-\log \frac{\Phi(e^{{\imunit}\vartheta})}{\Psi(e^{{\imunit}\vartheta})}-1\right\}d\vartheta$$
is the classical Itakura-Saito distance of maximum likelihood estimation for speech processing  \cite{DISTORTION,BASSEVILLE_DISTANCEMEASURES}.}
We { now formulate} the following {\em Spectrum Approximation Problem:}
\begin{probl}
\label{problem:formalization}
Let $\Psi\in \mathcal{S}_+^{m \times m}$, $G(z)$ as in \eqref{eq:filter} and $\Sigma=\Sigma^\top>0$. Find $\Phi^\circ$ that solves:
$$
 \rm{minimize} \; {d_{RER}}(\Phi,\Psi)\ \ \rm{over} \; \left\lbrace \Phi \in \Specspacepiu | \int G\Phi G^* = \Sigma \right\rbrace.
$$
\end{probl}

\begin{remark}
Notice that we could also minimize the distance index \eqref{RERDIST} with respect to the second argument. Indeed, this choice is meaningful in some approximation problems related to minimum prediction error and model reduction, see \cite{LPBook}. In our approximation problem (\ref{problem:formalization}), however , it is possible to prove that such a choice usually leads to a non rational approximant, even when the prior $\Psi$ is rational. Therefore, this approach is not suitable for our purposes.
\end{remark}

We first address the issue of  feasibility of Problem \ref{problem:formalization}, namely existence of $\Phi\in {\cal S}^{m\times m}_+(\T)$ satisfying (\ref{CONST})
where $G$ is the transfer function of the bank of filters (\ref{statefilter}) and $\Sigma$ is the steady-state covariance of the output process.
%
%
To this aim we first introduce some notation: All through the paper,  $\symm_m\subset \R^{m\times m}$ denotes the $m(m+1)/2$-dimensional,   real  vector space of $m$-dimensional
symmetric matrices.
We denote by $\Ctpiu$  the set of continuous spectral densities of $m$-dimensional $\R^m$-valued processes defined on the unit circle $\mathbb{T}$. We indicate by $\Ct$ the linear space generated by $\Ctpiu$. Let $\Gamma \, :\, \Ct \rightarrow \symm_n$ be the linear operator  defined by
\begin{equation}\label{eq:gamma}
 \Gamma (\Phi):= \int G\Phi G^*.
\end{equation}
The following result { can be obtained along the same lines of  \cite{GEORGIOU_STATECOV} (see also \cite{RAMPONI_FERRANTE_PAVON_GLOBALLYCONVERGENT}).\footnote{In \cite{GEORGIOU_STATECOV} the general case was considered when $A\,\in\C^{n\times n}$, $B\,\in\C^{n\times m}$ and the process $y$ is complex-valued, too. In that case, it was proven that the Hermitian matrix $\Sigma\in\C^{n\times n}$ belongs to $\Range(\Gamma)$ if and only if there exists $H \in \mathbb{C}^{m\times n}$ solving the feasibility equation $\Sigma-A\Sigma A^\ast=BH + H^\ast B^\ast$.}}
\begin{theorem}
 \label{thm:feasibilitycond}
 Consider $\Sigma = \Sigma\trasp \,\in \mathbb{R}^{n\times n}$ and a filter defined as in \eqref{eq:filter}. 
Then:
\begin{enumerate}
\item $\Sigma$ is in $\Range(\Gamma)$ if and only if there exists $H \in \mathbb{R}^{m\times n}$ such that
\begin{equation}
\label{eq:feasibility}
 { \Sigma-A\Sigma A\trasp=BH + H\trasp B\trasp}.
\end{equation}
\item Let the $\Sigma\,\in\R^{n\times n}$ be positive definite. Then, there exists $H \, \in \mathbb{R}^{m\times n}$ that solves \eqref{eq:feasibility} if and only if there exists $\Phi \, \in \Ctpiu$ such that $\Gamma(\Phi)=\Sigma$.
\end{enumerate}
\end{theorem}
 From now on we assume feasibility of Problem \ref{problem:formalization}.
In view of the previous result, this is equivalent to the fact that  Equation (\ref{eq:feasibility}) admits { a solution $\bar{H}$}.
Moreover, to simplify the exposition, we assume that $\Sigma=I$.
This can be done without loss of generality. In fact, if $\Sigma\neq I$, it suffices to replace $G$ with $G':=\Sigma^{-1/2}G$ and $(A,B)$ with
$(A'=\Sigma^{-1/2}A\Sigma^{1/2},B'=\Sigma^{-1/2}B)$ to obtain an equivalent problem
where $\Sigma= I$.
We  now proceed to solve Problem \ref{problem:formalization}. Since
\begin{equation}\label{eq:constterm}
\frac{1}{4\pi}\int_{-\pi}^\pi\left\{-\tr\Psi^{-1}(e^{{\imunit}\vartheta})\Psi(e^{{\imunit}\vartheta})\right\}d\vartheta=-\frac{m}{2}
\end{equation}
the left-hand side of (\ref{eq:constterm})
plays no role in the optimization. It can, therefore, be neglected together with a $\frac{1}{2}$  multiplying the integral. Thus, Problem  \ref{problem:formalization} is equivalent to minimizing, over  $\Specspacepiu$,
$2d_{RER}(\Phi,\Psi)+m=\int\left\{\log\det \left(\Phi^{-1}\Psi\right)+\tr\left(\Psi^{-1}\Phi\right)\right\},$
subject to (\ref{CONST}).  Recall that  the inner product in $\symm_n$ is defined by $\langle M,N\rangle=\tr[MN]$. We can then consider  the Lagrangian
\begin{eqnarray}\nonumber
L_\Psi(\Phi, \Lambda) &=&2d_{RER}(\Phi,\Psi)+m+\langle\Lambda,\int G\Phi G^* -\Sigma\rangle\\&=& \int \left[\log\frac{\det(\Psi)}{\det(\Phi)}+\tr(\Psi ^{-1}\Phi) + \tr(\Lambda G\Phi G^*)\right] -\tr\Lambda,\label{eq:lag}
\end{eqnarray}
where the Lagrange parameter $\Lambda \in \symm_n$ and we have used the assumption $\Sigma=I$.
Notice that each $\Lambda \in \symm_n$ can be uniquely  decomposed as $\Lambda = \Lambda_\Gamma + \Lambda_\perp$, where $\Lambda_\Gamma \in \Range{(\Gamma)}$ and $\Lambda_\perp \in \left(\Range{(\Gamma)}\right)^\perp$.
It can be proven \cite[Section III]{RAMPONI_FERRANTE_PAVON_GLOBALLYCONVERGENT} that, $\forall \, \Lambda_\perp \in {\left(\Range{(\Gamma)}\right)}^\perp$,
$G^*(e^{{\imunit}\vartheta}) \Lambda_\perp G(e^{{\imunit}\vartheta}) \equiv 0$. Moreover, $\tr{\left[\Lambda_\perp\right]}=\left\langle\Lambda_\perp,I\right\rangle=0$, because $I\, \in \Range{(\Gamma)}$ in view of the feasibility assumption.
Hence, a term $\Lambda_\perp \, \in \left(\Range{(\Gamma)}\right)^\perp$ gives no contribution to the Lagrangian \eqref{eq:lag}.  We therefore assume from now on that the Lagrange parameter $\Lambda$ belongs to $\Range(\Gamma)$.

For $\Lambda$ fixed, we consider now the \emph{unconstrained} minimization of the functional \eqref{eq:lag} with respect to $\Phi$.
Observe that   $L_\Psi(\cdot, \Lambda)$ in \eqref{eq:lag} is strictly convex on $ \mathcal{S}_+^{m \times m}$. We impose that the first variation  be zero in each direction $\delta\Phi\in L_2^{m\times m}$.
Recalling that, for a positive definite matrix $X$, the directional derivative of
$\log\det(X)$ in direction $\delta X$ is given by
\begin{equation}\label{eq:deltalog}
\delta\log\det(X;\delta X)=\tr(X^{-1}\delta X),
\end{equation}
 we get:
\begin{equation}\label{eq:deltaLag}
 \delta L(\Phi, \Lambda; \delta \Phi) = \int \left[ -\tr(\Phi^{-1}\delta\Phi)+\tr(\Psi^{-1}\delta\Phi)+\tr (G^*\Lambda G\delta\Phi)\right]= \int \langle-\Phi^{-1}+\Psi^{-1}+ G^*\Lambda G,\delta\Phi\rangle.
\end{equation}
Since $\left[-\Phi^{-1}+\Psi^{-1}+ G^*\Lambda G\right]\in L_2^{m\times m}$, \eqref{eq:deltaLag} is zero $\forall \, \delta\Phi\in L_2^{m\times m}$ if and only if
\begin{equation}\label{eq:phiOpt}
 \Phi= { \Phi^\circ(\Lambda)} := {\left[\Psi^{-1}+G^* \Lambda G\right]}^{-1}.
\end{equation}
Let $W_\Psi$ be the stable and minimum phase spectral factor of $\Psi$,\footnote{  Since $\Psi\in \Specspacepiu$, $W_\Psi$ exists. It is unique up to multiplication on the right by a constant orthogonal matrix.}  and   $G_1(e^{{\imunit}\vartheta})$ be defined by
\begin{equation}\label{eq:G1}
 G_1(e^{{\imunit}\vartheta}):= G(e^{{\imunit}\vartheta})W_\Psi(e^{{\imunit}\vartheta}).
\end{equation}
It will be later interesting to consider also the alternative form of (\ref{eq:phiOpt})
\eq\label{eq:alternative}
{\Phi^\circ(\Lambda)} = W_\Psi{(I+G_1 ^* \Lambda G_1)}^{-1}W_\Psi ^*.
\eeq
It is important to point out that  \eqref{eq:phiOpt} yields an upper bound on the  McMillan degree $\deg [\Phi^\circ]$ of the optimal approximant $\Phi^\circ$. Indeed, {it follows from  \eqref{eq:phiOpt} that $\deg [\Phi^\circ]\leq \deg [\Psi] + 2n$, where $n$ is the McMillan degree} of $G(z)$. This result represents a significant improvement in the frame of multivariable spectral estimation, {in which the best so far available upper bound on the McMillan degree (which can be regarded as a measure of complexity)  of the solution  was $\deg [\Psi]+4n$ (see \cite{FERRANTE_PAVON_RAMPONI_HELLINGERVSKULLBACK}).}

Since $\Phi^\circ$ is required to be a bounded spectral density, we need, as indicated by \eqref{eq:alternative}, to restrict the Lagrange multiplier $\Lambda$ to the subset $\mathcal{L}_+$, where
\begin{equation}\label{eq:lpiu}
 \mathcal{L}_+:=\left\lbrace\Lambda\in\,{\symm_n}\, |\, I+G_1^*\Lambda G_1 > 0 \, \, \text{ a.e. on}\; \mathbb{T}\right\rbrace.
\end{equation}
In conclusion, the natural set for the Lagrangian multiplier $\Lambda$ is
\begin{equation}\label{eq:Lpiugamma}
 \mathcal{L}_+ ^\Gamma := \mathcal{L}_+ \cap \Range{(\Gamma)}.
\end{equation}
To sum up, the main result is that for each $\Lambda\in\mathcal{L}_+ ^\Gamma$ there exists a unique $\Phi ^\circ\in \, \mathcal{S}_+ ^{m\times m}$ that minimizes the Lagrangian functional. It has the form \eqref{eq:phiOpt}.
If we produce a $\Lambda^\circ$ s.t.  $\Phi^\circ(\Lambda^\circ)$ satisfies  constraint (\ref{CONST}), then such a  $\Phi^\circ(\Lambda^\circ)$ is the solution of Problem \ref{problem:formalization}.
Existence of such a $\Lambda^\circ$ turns out to be a most delicate issue. To address this problem, we resort to duality.

\section{The dual problem}\label{dualproblem}
{
Consider
$$\inf _\Phi L(\Phi,\Lambda)=L(\Phi ^\circ, \Lambda)= \int \log \det (I+G_1 ^*\Lambda G_1) +n - \tr\Lambda.$$
Instead of maximizing this expression, we will equivalently minimize  the following functional hereafter referred to as the \emph{dual functional}:
\begin{equation}\label{eq:J}
J_{\Psi}(\Lambda) := -L(\Phi ^\circ(\Lambda),\Lambda)+n=\int \left[\tr\Lambda - \log\det (I+G_1 ^*\Lambda G_1) \right].
\end{equation}
{Recall that given a matrix $A=A^\ast>0$,
we have $\int \log\det{A} = \int \tr \log{A}$. Hence, we can  express the dual functional also as
$J_{\Psi}(\Lambda) = \int \tr \left[\Lambda-\log(I+G_1 ^*\Lambda G_1)\right].$}
Given $\delta \Lambda\in \symm_n$, by means of \eqref{eq:deltalog} we can evaluate its first variation:
\begin{equation}\label {eq:der1}
 \delta J_\Psi(\Lambda;\delta \Lambda)=\nabla J_{\Psi,\Lambda}(\delta \Lambda)= \int \left\lbrace \tr\left[\delta\Lambda\right] - \tr \left[{(I+G_1 ^*\Lambda G_1)}^{-1}G_1 ^*\delta \Lambda G_1\right]\right\rbrace.
\end{equation}

{ The results of this section show that there exists a unique $\Lambda^\circ \in \LL $ minimizing
$J_{\Psi}(\Lambda)$ in (\ref{eq:J}). Such a $\Lambda^\circ$ annihilates the directional derivative (\ref{eq:der1}) in any direction $\delta\Lambda\in
\symm_n$, namely
\begin{equation}
\langle I-\int G_1(I+G_1 ^*\Lambda^\circ G_1)^{-1}G_1 ^*,\delta\Lambda\rangle=0\ \ \
\forall \delta\Lambda\in
\symm_n,
\end{equation}
or, equivalently,
\begin{equation}
I=\int G_1(I+G_1 ^*\Lambda^\circ G_1)^{-1}G_1 ^*=\int G\Phi^\circ(\Lambda^\circ) G^*.
\end{equation}
This means that the corresponding spectral density $\Phi^\circ:=\Phi(\Lambda^\circ)= {\left[\Psi^{-1}+G^* \Lambda^\circ G\right]}^{-1}$, satisfies  constraint (\ref{CONST}) { (recall that we set $\Sigma =I$)} and is therefore the unique solution of  Problem \ref{problem:formalization}.}

{ Uniqueness of the minimizing $\Lambda_0\in \mathcal{L}_+ ^\Gamma$  is an obvious consequence of the following result.}
\begin{theorem}\label{theorem:contConv}
 The dual functional $J_\Psi(\Lambda)$ belongs to  $\mathcal{C}^2(\mathcal{L}_+ ^ \Gamma)$ and is \emph{strictly} convex on $\mathcal{L}_+ ^ \Gamma$.
\end{theorem}
\begin{IEEEproof}
Consider a sequence $M_n \,\in \, \Range{(\Gamma)}$, such that $M_n\rightarrow \, 0$, and define, for $N\in\symm_n$,  $Q_N(z)=I+G_1^*(z)NG(z)$.  By Lemma 5.2 in \cite{RAMPONI_FERRANTE_PAVON_GLOBALLYCONVERGENT},  ${Q}^{-1} _{\Lambda + M_n}$ converges uniformly to ${Q^{-1} _{\Lambda}}$, so that it is  bounded above. Hence, applying the bounded convergence theorem, we get
$$
\lim _{n \rightarrow \infty} \int \tr \left[{Q}^{-1} _{\Lambda + M_n} G_1 ^*\delta \Lambda G_1\right] = \int \tr \left[{Q}^{-1} _{\Lambda} G_1 ^*\delta \Lambda G_1\right],
$$
so that $J_\Psi(\Lambda)$ belongs to  $\mathcal{C}^1(\mathcal{L}_+^\Gamma)$.
Consider now the second variation. Let us denote the matrix inversion operator by $R:M\mapsto M^{-1}$ and  recall that its first derivative in direction $\delta M$ is given by $\delta R \left(M,\delta M\right)=-M^{-1}\delta M M^{-1}$.   Then, for  $\delta \Lambda_1$ and $\delta\Lambda_2$ in $\symm_n$, we have
\begin{equation}\label {eq:der2}
 \delta ^2 J_\Psi(\Lambda;\delta \Lambda_1, \delta \Lambda_2)= \int  \tr \left[{(I+G_1 ^* \Lambda G_1 )}^{-1}G_1 ^*\delta \Lambda _2 G_1{(I+G_1 ^* \Lambda G_1 )}^{-1}G_1 ^*\delta \Lambda _1 G_1\right],
\end{equation}
so that $J_\Psi(\Lambda)$ is $\mathcal{C}^2(\mathcal{L}_+^\Gamma)$. The bilinear form $H_\Lambda(\cdot, \cdot) := \delta^2\Jpsi(\Lambda; \cdot, \cdot)$
is the {\em Hessian} of $\Jpsi$ at $\Lambda$. For $\delta\Lambda\in\Range(\Gamma)$, which implies that $(\Lambda+\varepsilon\delta\Lambda)\in\mathcal{L}_+ ^ \Gamma$ for sufficiently small $\varepsilon$, consider $H_\Lambda(\delta\Lambda,\delta\Lambda)= \delta ^2 J_\Psi(\Lambda;\delta \Lambda, \delta \Lambda)$.  We get
\begin{equation}\label{eq:strictConv}
\begin{split}
H_\Lambda(\delta\Lambda,\delta\Lambda)
                                     &= \int  \tr \left[{(I+G_1 ^* \Lambda G_1 )}^{-1}G_1 ^*\delta \Lambda G_1{(I+G_1 ^* \Lambda G_1 )}^{-1}G_1 ^*\delta \Lambda G_1\right]\\
                                     &=\int \tr \left[{Q^{-\frac{1}{2}}_\Lambda}G_1 ^*\delta \Lambda G_1{Q}^{-1} _{\Lambda}G_1 ^*\delta \Lambda G_1 {Q^{-\frac{1}{2}}_\Lambda}\right]
\end{split}
\end{equation}
which vanishes if and only if the integrand is identically zero. Moreover $G_1 ^*\delta\Lambda G_1 = W_\Psi ^*G^*\delta\Lambda G W_\Psi$ is identically zero on $\T$ if and only if $\delta\Lambda \,\in\, \Range(\Gamma)^\perp$.
On the other hand we have assumed $\delta\Lambda\in \Range(\Gamma)$, so that the integrand is identically zero if and only if $\delta\Lambda =0$. In conclusion, the Hessian is positive-definite and  the dual functional is strictly convex on  $\mathcal{L}_+ ^ \Gamma$.
\end{IEEEproof}

The next and most delicate step is to prove that, although the set $\mathcal{L}_+^\Gamma$ is open and unbounded,  a $\Lambda^\circ$ minimizing $J_\Psi$ over $\mathcal{L}_+ ^\Gamma$ \emph{does} exist.
To this aim, first we prove that the function $J_\Psi(\Lambda)$ is inf-compact,  i.e. $\forall \, \alpha \in \mathbb{R}$, the set $
\left\lbrace\Lambda \in \mathcal{L}_+ ^\Gamma \, | \, J_\Psi(\Lambda)\leq \alpha\right\rbrace
$ is compact.
To establish this fact, define $\CL$ to be the closure of $\mathcal{L}_+^\Gamma$, i.e. the set
 $$
 \CL = \left\lbrace\Lambda=\Lambda\trasp \in \mathbb{R}^{n \times n} \, | \, \Lambda \in \Range(\Gamma),\, I+G_1^*\Lambda G_1\geq 0, \, \forall e^{{\imunit}\vartheta} \in \mathbb{T} \right\rbrace.
 $$
Given that, for $\Lambda$ belonging to the boundary $\partial \mathcal{L}_+^\Gamma$, the Hermitian matrix $I+G_1^*\Lambda G_1$ is singular, in at least one point of $\T$, it is useful to introduce the following sequence of functions on $\CL$:
\begin{equation}
 J_\Psi ^n(\Lambda) = \int \tr\left[\Lambda - \log\left(I+G_1^*\Lambda G_1 + \frac{1}{n}I\right)\right],\quad n\ge 1.
\end{equation}
Recall that a real-valued function $f$ is said to be lower semicontinuous at $x_0$ if, $\forall \, \varepsilon >0$, there exists a neighborhood $U$ of $x_0$ such that, $\forall \, x \, \in U$, $f(x)\geq f(x_0)-\varepsilon$.
Recall also that, given $f: \mathbb{R}^{n\times n}\rightarrow \mathbb{R}$, its \emph{epigraph} $\epi(f)$ is defined by
$$\epi(f) := \left\lbrace (x,a) \, \in\, \mathbb{R}^{n\times n}\times \mathbb{R} \,|\,  a\geq f(x)\right\rbrace.$$
Moreover, $f$ is a lower semicontinuous (convex) function if and only if its epigraph is closed (convex), see e.g.  \cite{ROCKAFELLAR_CONVEX}.
The following Lemmata allow to conclude that $J_\Psi(\Lambda)$ is inf-compact over $\CL$.
\begin{lemm}\label{lemma:lsc} 
The pointwise limit $J_\Psi ^\infty(\Lambda)$, defined as $J_\Psi^\infty(\Lambda):=\lim_{n\rightarrow \infty}J_\Psi^n(\Lambda)$, exists and is a lower semicontinuous and convex function defined over $\CL$, with values in the extended reals.
\end{lemm}
\begin{IEEEproof}
The additive term $\frac{1}{n}I$ ensures that, for each $n$, $J_\Psi^n(\Lambda)$ is a continuous and convex function of $\Lambda$ on the closed set $\CL$.
From the properties of $J_\Psi^n(\Lambda)$, it follows that $\epi(J_\Psi^n(\Lambda))$ is a closed and convex subset of $\mathbb{R}^{n\times n}\times \mathbb{R}$. In addition, the pointwise sequence is monotonically increasing, since $J_\Psi^n(\Lambda)<J_\Psi^{n+1}(\Lambda)$. Therefore, it converges to $J_\Psi^\infty(\Lambda):=\sup_n{J_\Psi^n(\Lambda)}$. Since the intersection of closed sets is closed and the intersection of convex sets is convex, $\epi J_\Psi^\infty(\Lambda)=\cap_n \epi J_\Psi^n(\Lambda)$ is closed and convex. As a consequence, $J_\Psi^\infty(\Lambda)$ is lower semicontinuous and convex. \end{IEEEproof}
\begin{lemm}\label{lemm:tracciaLambda}
Assume that the feasibility condition \eqref{eq:feasibility} holds. { Given $\Lambda \in \mathcal{L}_+^\Gamma$,} there exist two real constants $\mu>0$ and $\alpha$ such that:
\begin{equation}
\tr\left[\Lambda\right] \geq \mu \tr \left[\int (G_1^*\Lambda G_1+I)\right] +\alpha.
\end{equation}
\end{lemm}
\begin{IEEEproof}
Since $\Sigma=I$, by feasibility, there exists $\Phi_I\in\Specspacepiu$ such that $\int G\Phi_I G^*=I$. Thus,
\begin {equation}\label{eq:trLambda}
\begin{split}
  \tr\left[\Lambda\right] &= \tr\left[\int G\Phi_IG^* \Lambda\right] = \tr\left[ \int  G^*\Lambda G \Phi_I\right]\\
  &= \tr\left[ \int W_\Psi^{*} G^*\Lambda G W_\Psi W_\Psi^{-1} \Phi_I W_\Psi^{-*}\right] = \tr\left[ \int  G_1^*\Lambda G_1 \Xi \right],
  \end{split}
\end {equation}
where the cyclic property of the trace was employed and the auxiliary spectral density $\Xi:= W_\Psi ^{-1}\Phi_I W_\Psi^{-*}$ has been defined.
By defining $\alpha:=- \tr\left[\int \Xi \right]$, it follows that
\begin{equation}\label{eq:trLambdaPart2}
\begin{split}
\tr\left[\Lambda\right]
            &= \tr\left[\int (G_1^*\Lambda G_1+I)\Xi\right] - \tr\left[\int \Xi \right]= \tr\left[\int (G_1^*\Lambda G_1+I)\Xi\right]+\alpha.
 \end{split}
\end {equation}
Let $\Delta$ be such that $(G_1^*\Lambda G_1+I)=\Delta^*\Delta$
 (recall that we are assuming $\Lambda\in\LL$ so that $G_1^*\Lambda G_1+I$ is positive definite on $\T$ and admits a right spectral factor $\Delta$)
so that $\tr\left[ (G_1^*\Lambda G_1+I)\Xi\right] = \tr\left[ \Delta\Xi\Delta^*\right].$
Given that $\Xi =W_\Psi ^{-1}\Phi_I W_\Psi^{-*}$ is a coercive spectrum, because both $\Phi_I$ and $\Psi$ belong to $\Specspacepiu$, there exists $\mu>0 \, \text{  s.t.    } \, \Xi(e^{{\imunit}\vartheta}) \geq \mu I, \, \forall \, e^{{\imunit}\vartheta} \in \mathbb{T}$. Recalling that the trace and the integral are monotonic functionals, it is possible to conclude that
\begin{equation}\label{eq:trLambdaPart3}
  \tr\left[\Lambda\right]= \tr\left[\int (G_1^*\Lambda G_1+I)\Xi\right] +\alpha \geq \mu \tr \left[\int (G_1^*\Lambda G_1+I)\right] +\alpha.\end{equation}
\end{IEEEproof}

\begin{lemm}\label{lemma:bd}
Define $\mathcal{B} := \left\lbrace \Lambda \, \in \, \partial\LL \, | \,\det{\left(G_1^*\Lambda G_1 + I\right)} = 0, \, \forall e^{{\imunit}\vartheta} \,\in \mathbb{T}\right\rbrace$ and consider its complement set $\mathcal{B}^c := \left\lbrace\Lambda \in \partial\LL \, |\, \Lambda \notin \mathcal{B}\right\rbrace $. Then,  under feasibility assumption:
\begin{enumerate}
\item $J_\Psi^\infty(\Lambda)$ is { bounded from below} on $\CL$;
\item $J_\Psi^\infty(\Lambda)=J_\Psi(\Lambda)$ on $\LL$;
\item $J_\Psi^\infty(\Lambda)$ is finite over $\mathcal{B}^c$.
\end{enumerate}
\end{lemm}
The proof can be found in the Appendix.

\begin{lemm}\label{lemma:infJ}
 If the feasibility hypothesis holds, then, for $\Lambda\in\LL$,
\begin{equation}
 \lim_{\|\Lambda\|\rightarrow + \infty} J_\Psi(\Lambda)= +\infty.
\end{equation}
\end{lemm}
 See the Appendix for the proof.
Then, by Weierstrass' Theorem we can conclude that there exists a minimum point $\Lambda^\circ\in \CL$. More can be proven:
\begin{theorem}\label{theorem:infC}
If the feasibility condition (\ref{eq:feasibility}) holds, the problem of minimizing $J_\Psi(\Lambda)$ over $\LL$ admits a unique solution $\Lambda^\circ\, \in\, \LL$.
\end{theorem}
\begin{IEEEproof}
Since $J_\Psi(\Lambda)$ is inf-compact over $\CL$, it admits a minimum point $\Lambda^\circ$ there. Obviously, $\Lambda^\circ\notin \mathcal{B}$, since $J_\Psi(\Lambda) = +\infty$ on $ \mathcal{B}$ (Lemma \ref{lemma:bd}). Suppose $\Lambda^\circ \, \in \, \mathcal{B}^c$. By Lemma \ref{lemma:bd} again,  it follows that $J_\Psi(\Lambda^\circ)$ is finite. By convexity of $\CL$, $\forall \, \varepsilon \, \in \, \left[0,1\right]$, $\Lambda^\circ + \varepsilon(I-\Lambda^\circ) \, \in\, \CL$, since the feasibility condition \eqref{eq:feasibility} ensures that $I\in \LL$. The one-sided directional derivative is
 \begin{equation}
  \begin{split}
   \delta J_{\Psi_+}(\Lambda^\circ; I-\Lambda^\circ)&=\lim_{\varepsilon \searrow 0} \left[\frac{J_\Psi\left(\Lambda^\circ+\varepsilon\left(I-\Lambda^\circ\right)\right)-J_\Psi(\Lambda^\circ)}{\varepsilon}\right]\\
 &= \tr\left[I-\Lambda^\circ\right] - \int \tr \left[{\left(I+G_1^*\Lambda^\circ G_1\right)}^{-1}G_1^*\left(I-\Lambda^\circ\right)G_1\right]\\
 &= \tr\left[I-\Lambda^\circ\right] - \int \tr \left[{\left( I+ G_1^*\Lambda^\circ G_1\right)}^{-1}\left(G_1^*G_1-G_1^*\Lambda^\circ G_1 +I -I\right)\right]\\
 &= \tr\left[I-\Lambda^\circ\right] - \int \tr \left[{\left( I+ G_1^*\Lambda^\circ G_1\right)}^{-1}\left(I+G_1^*G_1\right) -I\right]\\
 &=-\infty.
  \end{split}
 \end{equation}
The last equality holds because for each ${\Lambda}\in \mathcal{B}^c$, the matrix $I+ G_1^*{\Lambda}G_1$ is singular and $I+ G_1^*G_1>0$ on $\T$.
As a consequence, the minimum point cannot belong to $\BL$. Thus, $\Lambda^\circ \, \in \LL$.
\end{IEEEproof}

Finally, we are left with the problem of developing an efficient numerical algorithm to compute the optimal solution $\Lambda^\circ$.
}

\section{Efficient implementation of a matricial Newton-like algorithm}
\label{nnsec6}

In order to compute the minimizer of the dual functional $J_{\Psi}(\Lambda)$, a matricial Newton-type algorithm is proposed. Here are the main steps: { (i) the starting point for the minimizing sequence $\left\lbrace\Lambda_i\right\rbrace_{i\in \mathbb{N}}$ is $\Lambda_0=0$, (ii) at each step we compute the Newton \emph{search direction} $\Delta \Lambda _i$, (iii) we compute the Newton \emph{step length} $t_i^k$.}

\subsection{Search Direction}
{ Even though the problem is finite dimensional, the computation of the search direction is
rather delicate}  because a matricial expression of the Hessian and the gradient allowing to compute the search direction $\Delta x$ as $\Delta x = - {H_x}^{-1}\nabla f_x$ is not available.
In order to compute $\Delta \Lambda _i$, given $\Lambda_i\in\mathcal{L}_+ ^\Gamma$, one has to solve,
for the unknown $\Delta\Lambda_i$, the equation
$
H_{\Lambda _i}(\Delta \Lambda _i, \cdot)=-\nabla J_{\Psi,{\Lambda _i}}(\cdot),
$
which can be explicitly written as:
\begin{equation*}
 \int G_1{(I+G_1 ^* \Lambda _i G_1)}^{-1}G_1 ^* \Delta \Lambda _i G_1{(I+G_1 ^* \Lambda _i G_1)}^{-1}G_1 ^*=\int G_1{(I+G_1 ^* \Lambda _i G_1)}^{-1}G_1 ^*-I.
\end{equation*}

To this aim, consider a basis of $\Range(\Gamma)$. It can be readily obtained, by recalling that $\Sigma _k \in \, \Range(\Gamma) $ if and only if $ \exists \,H_k \in \mathbb{R}^{m \times n}$ s.t. $\Sigma _k - A\Sigma _k A\trasp = BH_k + {H_k}\trasp B\trasp $.
Therefore, considering a basis $\left\lbrace H_1,\dots, H_L\right\rbrace$ for $\mathbb{R}^{m \times n}$, a set of generators $\left\lbrace \Sigma_1',\dots, \Sigma_L'\right\rbrace$ can be found by solving $L$ Lyapunov equations.
After that a basis $\left\lbrace \Sigma_1',\dots, \Sigma_N'\right\rbrace$ can be easily computed.\footnote{Indeed, following the lines detailed in  \cite{FERRANTE_PAVON_ZORZI_ENHANCEMENT}, it is possible to obtain {\em directly} a basis of $\Range(\Gamma)$  by solving only $N$ Lyapunov equations.}
Since $I\, \in\,\Rgamma$, we can add to each $\Sigma_i$ the matrix $\alpha_i I$, and, for suitable (large) $\alpha_i$, get a basis $\left\lbrace \Sigma_1,\dots, \Sigma_N\right\rbrace$ of $\Range(\Gamma)$ made of positive definite matrices.
The search direction can now be computed by applying the following procedure:
\begin{enumerate}
\item Compute
      \begin{equation}\label{eq:Y}
      Y=\int G_1{(I+G_1 ^* \Lambda _i G_1)}^{-1}G_1 ^*-I
      \end{equation}
\item For each generator $\Sigma_k$, compute
      \begin{equation}\label{eq:Yk}
      Y_k=\int G_1{(I+G_1 ^* \Lambda _i G_1)}^{-1}G_1 ^* \Sigma _k G_1{(I+G_1 ^* \Lambda _i G_1)}^{-1}G_1 ^*
      \end{equation}
\item Find $\left\lbrace\alpha _k\right\rbrace$ s.t. $Y=\sum _k \alpha _k Y_k$;
\item Set $\Delta\Lambda _i = \sum_k \alpha _k \Sigma _k$.
\end{enumerate}
The most challenging step is to compute $Y$ and $Y_k$. A sensible approach is to employ spectral factorization techniques in order to compute the integrals, along the same lines described in \cite[Section VI]{RAMPONI_FERRANTE_PAVON_GLOBALLYCONVERGENT}.
Indeed, the integrand that appears in equation \eqref{eq:Y} is a coercive spectral density and the same holds for the integrand in \eqref{eq:Yk}, since we have chosen  the generators $\Sigma_i$ to be positive definite. As a consequence, the integral may be computed by means of numerically robust spectral factorization techniques.
For the computation of $Y$, let us focus on $Q_{\Lambda_i} (z)= I+G_1 ^* (z) \Lambda_i G_1(z)$.
Assume that a realization of the stable minimum phase spectral factor $W_\Psi(z)$ is given (or has been computed from $\Psi$). Then, we can easily obtain a state-space realization $G_1(z)=C_1(zI-A_1)^{-1}B_1$ of $G_1$.
Since $\Lambda_i\in\mathcal{L}_+ ^\Gamma$, $Q_{\Lambda_i} (z)$ is positive definite on $\mathbb{T}$, so that  the
following ARE admits a positive definite stabilizing solution $P=P\trasp>0$
(see, e.g. Lemma 6.4 in \cite{RAMPONI_FERRANTE_PAVON_GLOBALLYCONVERGENT}):
 \begin{equation}\label{eq:are}
   P=A_1\trasp PA_1 - A_1\trasp P B_1{(B_1\trasp P B_1 + I)}^{-1}B_1 \trasp PA_1+C_1\trasp \Lambda_i C_1.
 \end{equation}
Moreover,   $Q_{\Lambda_i}(z)$ can be factorized as $Q_{\Lambda_i}(z)=\Delta_{\Lambda_i} ^* (z)\Delta_{\Lambda_i}(z)$, where $\Delta_{\Lambda_i} (z)$ can be explicitly written in term of the stabilizing solution $P$:
\begin{equation}\label{eq:Delta}
 \Delta_{\Lambda_i} (z) =  (B_1 \trasp PB_1 + I)^{-\frac{1}{2}}B_1  \trasp PA_1  (zI-A_1)^{-1}B_1 + {(B_1 \trasp PB_1 + I)}^{\frac{1}{2}}.
\end{equation}
It is now easy to compute a state space realization of
$\Delta_{\Lambda_i}^{-1}$ and then of the stable filter
 $W_Y:=G_1\Delta_{\Lambda_i}^{-1}=C_1 (zI-Z_1)^{-1}B_1 (B_1 \trasp PB_1 + I)^{-\frac{1}{2}}$, with $Z_1:=A_1 - B_1{(B_1 \trasp PB_1 + I)}^{-1}B_1 \trasp PA_1$ being the closed-loop matrix.
The computation of (\ref{eq:Y}) is now immediate. In fact,
\begin{equation}\label{eq:factY}
Y +I = \int G_1{(I+G_1 ^* \Lambda _i G_1)}^{-1}G_1 ^* = \int G_1\Delta_{\Lambda_i} ^{-1}\Delta _{\Lambda_i} ^{-*}G_1 ^*= \int W_YW_Y ^*.
\end{equation}

The latter integral is thus the steady-state covariance of the output of the stable filter $W_Y$ driven by normalized white noise. It can be obtained by computing the unique solution of the Lyapunov equation $R-Z_1RZ_1 \trasp=B_1{(B_1 \trasp PB_1 + I)}^{-1}B_1\trasp$ and setting $Y +I =C_1RC_1\trasp$, so that
\begin{equation}Y =C_1RC_1\trasp-I .
\end{equation}
A similar procedure may be employed to compute also the matrices $Y_k$.

\subsection{Step length}
 The backtracking line search is implemented by halving the step
$t_i$  until both the following conditions are satisfied:
\eqn
&\Lambda _i + t_i ^k \Delta\Lambda _i \in \mathcal{L}_+ ^\Gamma; \label{eq:bt2}\\
 &J_\Psi (\Lambda _i + t_i ^k \Delta\Lambda _i)< J_\Psi (\Lambda _i) + \alpha t_i ^k \nabla J_{\Psi,\Lambda _i }\Delta\Lambda _i, \quad \text{where}\quad 0<\alpha<0.5.
\eeqn
The first condition can be easily evaluated by testing whether $Q_{\Lambda _i + t_i ^k\Delta\Lambda _i}$ admits a factorization of the kind introduced in the previous subsection or, equivalently, whether the corresponding ARE \eqref{eq:are} admits a solution $P=P\trasp>0$.

The only difficulty in checking the second condition is in computing
\begin{equation}
J_\Psi(\Lambda) = \tr \int \left[ \Lambda - \log(I + G_1 ^*\Lambda G_1)\right]=
\tr\Lambda -\int \log\det(I + G_1 ^*\Lambda G_1).
\end{equation}
The evaluation of the latter integral can be attained straightforwardly in the light of the fundamental result in statistical filtering (\ref{SK}).
In our case $Q(z)=Q_\Lambda(z)$ may be factorized as
$Q_\Lambda=\Delta^*\Delta$, where  $\Delta$ is a stable and minimum phase  filter for which a minimal realization can be computed as in the previous section (see eq. (\ref{eq:Delta})).
Since
$\log\det Q_\Lambda =\log\det\left[\Delta^*\Delta\right]=\log\det\left[\Delta\Delta^*\right]$,
$\det R$ is given by $\det[\Delta(\infty)\Delta^\ast(\infty)]$
which may be explicitly written in terms the solution $P$ of the corresponding ARE as  $\det[B_1 \trasp PB_1 + I]$.
Therefore,
$$
\int \log\det(I + G_1 ^*\Lambda G_1)=\log\det \left(B_1 \trasp PB_1 + I\right).
$$

\subsection{Convergence of the Proposed Algorithm}
A sufficient condition for global convergence of the algorithm is that the following requirements are  satisfied \cite[Chapter 9]{BOYD_VANDENBERGHE}:
\begin{enumerate}
\item $J_\Psi(\cdot)$ is twice continuously differentiable;
\item $\Lambda_0\in \mathcal{L}_+^\Gamma$ and the sublevel set $S:=\left\lbrace \Lambda \in \mathcal{L}_+^\Gamma| J_\Psi(\Lambda)\leq J_\Psi(\Lambda_0)\right\rbrace$ is closed;
\item $J_\Psi(\cdot)$ is \emph{strongly} convex, i.e. $\exists$ $m$ s.t. $H(J_\Psi)(\Lambda)>mI$, $\forall$ $\Lambda\in S$.
\item The Hessian is Lipschitz continuous in $S$, i.e. $\exists L$ such that:
      $$
      \begin{Vmatrix}H_{\Lambda_1}-H_{\Lambda_2}\end{Vmatrix}_2<L\begin{Vmatrix}\Lambda_2-\Lambda_1\end{Vmatrix}_2 \quad \forall \Lambda_1,\Lambda_2 \in S.
      $$
\end{enumerate}

In this case,  it is possible to prove not only that the algorithm converges, but also that, after a certain number of iterations, the backtracking line search always selects the full step (i.e. $t=1$). During the last stage the rate of convergence is quadratic, since there exists a constant $C$ such that $
\|\Lambda_{i+1}-\Lambda^\circ\|\leq C\|\Lambda_{i}-\Lambda^\circ \|^2.
$
Let us examine the requirements one by one.
The continuous differentiability of the dual function has already been proven in Section \ref{dualproblem}.
Theorem \ref{theorem:infC} states that the sublevel sets of { the dual function} $J_\Psi$ are compact, and hence closed (recall that, in a finite dimensional vector space, a set is compact if and only if it is closed and bounded).
Moreover, it is possible to conclude straightforwardly on \emph{strong} convexity and Lipschitz continuity of the Hessian.
Indeed, let us consider the sublevel set
$$
S=\left\lbrace\Lambda \in \mathcal{L}_+^\Gamma \, | \, J_\Psi(\Lambda)\leq J_\Psi(\Lambda _0)\right\rbrace.
$$
Notice that, assuming that $\Lambda_0$ is the starting point, the minimizing sequence computed by the Newton algorithm with backtracking line search is such that, $\forall \, k \geq 0,\ \Lambda_k\,\in\, S$.
The continuity of the Hessian over $\mathcal{L}_+^\Gamma$ has already been proven in Section \ref{dualproblem}. Moreover, since the map from a Hermitian matrix to its minimum eigenvalue is continuous (see Lemma 5.1 in \cite{RAMPONI_FERRANTE_PAVON_GLOBALLYCONVERGENT}), the map from $\Lambda \in \mathcal{L}_+^\Gamma$ to the minimum eigenvalue of $H_\Lambda(\delta\Lambda,\delta\Lambda)$ is continuous, being a composition of continuous maps. Since $S$ is compact, Weierstrass' Theorem holds. Therefore, there exists a minimum  $m$ in the set of eigenvalues of the Hessian $H_\Lambda(\delta\Lambda,\delta\Lambda),\,\forall\,\Lambda \in S$. Recall that the hypothesis of \emph{strict} convexity holds (as proven in Theorem \ref{theorem:contConv}). As a consequence, the Hessian $H_\Lambda$ is a positive definite matrix $\forall \, \Lambda \in \,S$, therefore $m>0$. In conclusion, there exists $m>0$ such that $H_\Lambda>mI, \, \forall \, \Lambda\, \in \,S$, i.e. $J_\Psi(\Lambda)$ is \emph!
 {strongly} convex.
Concerning the Lipschitz continuity of the Hessian of $J_\Psi(\Lambda)$, it is easy to see that $H_\Lambda$ is $\mathcal{C}^1(\mathcal{L}_+^\Gamma)$.
Indeed the third variation $\delta ^3 J_\Psi(\Lambda;\delta \Lambda_1, \delta \Lambda_2, \delta \Lambda_3)$ can be explicitly computed and its continuity can be proven along the same line developed in the proof of Theorem \ref{theorem:contConv} (the result can be extended, leading to the conclusion that $J_\Psi(\Lambda)$ is $\mathcal{C}^\infty(\mathcal{L}_+^\Gamma)$.
Continuous differentiability implies Lipschitz continuity on a compact set. Therefore, the Hessian is Lipschitz continuous on $S$.

In  conclusion, global convergence of the Newton algorithm is guaranteed, so that the proposed procedure is an effective computational tool to solve the spectral estimation Problem \ref{problem:formalization}.

\section{Simulation Results}\label{simul}
We now employ our results in a spectral estimation procedure, that may be outlined as follows.
\begin{enumerate}
\item We start from a finite sequence $\{y_1,\dots,y_N\}$, extracted from a realization of the zero-mean Gaussian process $y=\{y_k; k\in\mathbb{Z}\}$ with values in $\mathbb{R}^m$, whose spectrum is $\Phi(\ejth)$.
\item Design a filter $G(z)$, as described by equation \eqref{statefilter}.
\item Feed the filter with the data sequence $\{y_1,\dots,y_N\}$, collect the output data ${x_i}$ and compute a consistent estimate $\hat{\Sigma}$ of the covariance matrix.
\item In general, since the data length is finite, the estimate $\hat{\Sigma}$ does not satisfy the conditions stated in Theorem \ref{thm:feasibilitycond}. Our choice is to guarantee feasibility, by choosing a positive definite covariance matrix $\tilde{\Sigma}\in\Range{\Gamma}$ such that it is close to the starting estimate $\hat{\Sigma}$ in a suitable sense. Such an approach, introduced in \cite{FERRANTE_PAVON_ZORZI_ENHANCEMENT}, is briefly described in Remark \ref{remark:sigmastimata}.
\item Choose a prior spectral density $\Psi$.
\item Tackle Problem \ref{problem:formalization} by means of the proposed algorithm with the chosen $\Psi$ and $\Sigma=\tilde{\Sigma}$.
\end{enumerate}
\begin{remark}\label{remark:sigmastimata}
As previously observed, the covariance estimate $\hat{\Sigma}$ does not usually satisfy the feasibility requirements stated by Theorem \ref{thm:feasibilitycond}. In order to apply our method, we need to ensure feasibility. Thus, we need to compute an auxiliary positive definite covariance matrix $\tilde{\Sigma}$ which satisfies equation \eqref{CONST} and it is ``close'' to the estimate $\hat{\Sigma}$. To this purpose, an ancillary optimization problem is defined, so that the ``best'' approximant $\tilde{\Sigma}$ is chosen as
 \eq
\tilde{\Sigma} =  \underset{\Sigma\,\in\,\Range{\Gamma},\\ \Sigma>0}{\argmin} \D(p \|\hat{p}),
\eeq
where $p$ and $\hat{p}$ are zero-mean Gaussian densities with covariance $\Sigma$ and $\hat{\Sigma}$, respectively. The distance index is namely defined as follows:
\eq
\D (p\|\hat{p}) =\frac{1}{2} \left[\log\det{{\Sigma}^{-1}\hat{\Sigma}} + \tr{\hat{\Sigma}}^{-1}\Sigma -n\right].
\eeq
This problem can be solved efficiently by means of a matricial Newton algorithm. The reader is referred to \cite{FERRANTE_PAVON_ZORZI_ENHANCEMENT} for a complete exposition.
\end{remark}

Notice that our approach provides two degrees of freedom, since both the prior $\Psi$ and the filter $G(z)$ can be suitably selected.
Concerning the former, it can be chosen to be a coarse estimate $\hat{\Phi}$ of $\Phi$ obtained by means of a standard, ``simple'' estimation method. For instance, $\hat{\Phi}$ could be a low order ARMA model, computed through PEM methods. Such a choice gives a further insight into the meaning of the proposed procedure: Problem \ref{problem:formalization} consists in computing the bounded and coercive spectral density which is consistent with the interpolation condition \eqref{CONST} and is as close as possible to the initial estimate $\hat{\Phi}$ in the distance \eqref{relentropyconnection}.
Consider now the design of the filter. Recall that its role is to provide interpolation conditions for the approximant: In the light of this consideration, it is possible to reinterpret our approach as a generalization of classical problems such as Nevanlinna-Pick interpolation and the covariance extension problem, as explained in \cite[Section I]{BYRNES_GEORGIOU_LINDQUIST_THREE}.

\subsection{Scalar Case}
To begin with, we analyze the resolution capabilities of the proposed approach, by comparing its performances with those of the { original THREE method}.
In \cite{BYRNES_GEORGIOU_LINDQUIST_THREE}, it is explained how an adequate choice of the filterbank poles can improve the  estimate's resolution: a higher resolution can be attained by selecting poles in the proximity of the unit circle, with arguments in the range of frequency of interest.
\begin{figure}
\centering
\includegraphics[width=0.74\textwidth]{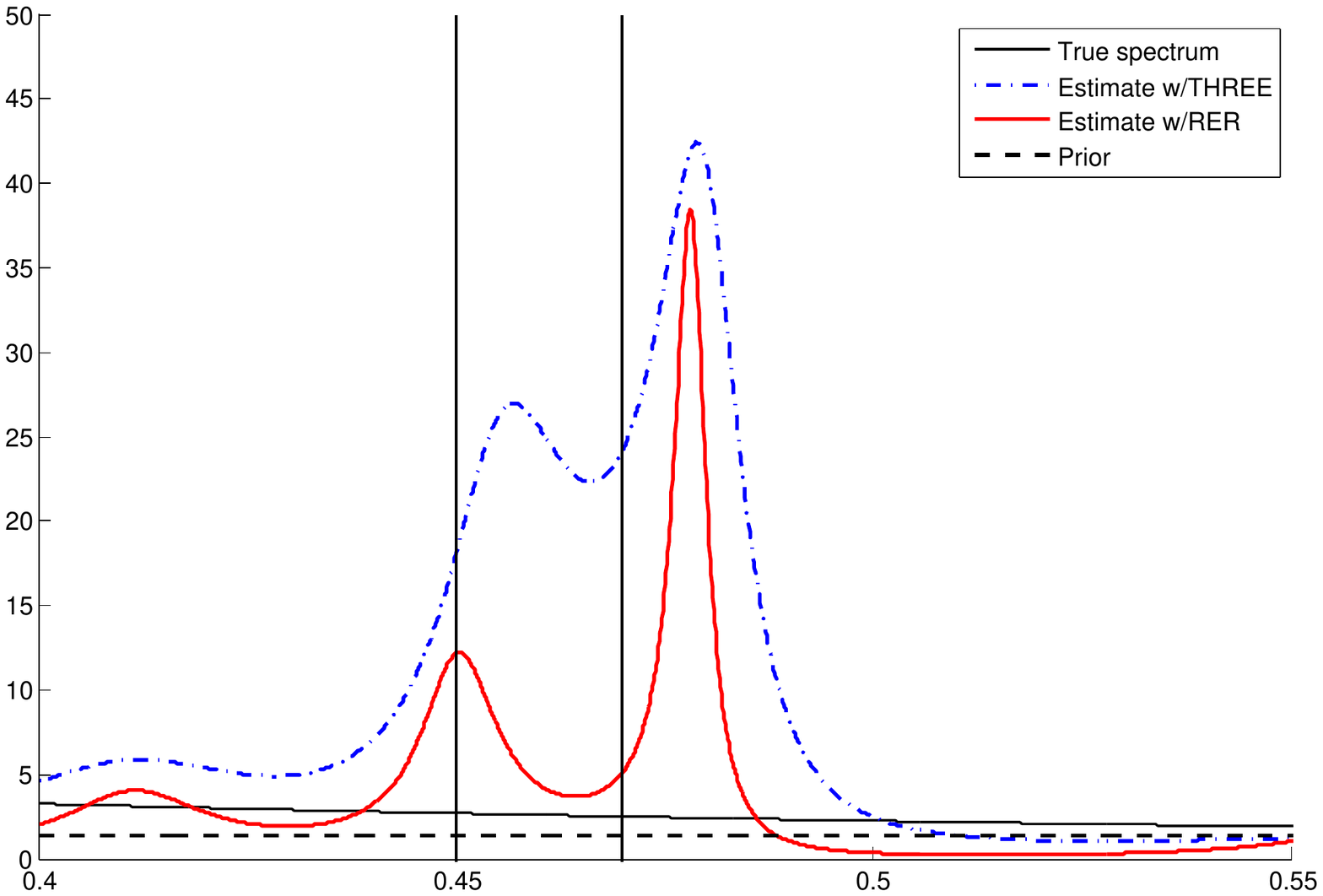}
\caption{Estimation of close spectral lines in colored noise ($\omega_1 = 0.45$ rad/s and $\omega_2 = 0.47$ rad/s). The chosen prior is the sample covariance of the data $\{y_k\}$. The  radius of the complex poles is equal to $0.95$. Both RER and THREE indicate the presence of the two lines.
}
\label{fig:testRigheVicine}
\end{figure}
In order to analyze such a feature, we deal with an instance of the classical problem of detecting spectral lines  in colored noise. The setting is the same described in \cite[Section IV.B]{BYRNES_GEORGIOU_LINDQUIST_THREE}.
The process of interest obeys to the following difference equation:
$$y(t)=0.5\sin(\omega_1t+\phi_1)+0.5\sin(\omega_2t+\phi_2)+z(t), \quad
z(t)=0.8z(t-1)+0.5\nu(t)+0.25\nu(t-1),
$$
where the variables $\phi_1$, $\phi_2$ and $\nu(t)$ are Gaussian, independent, with zero-mean and unit variance.
Matrix $B$ is a column of ones. Matrix $A$ was chosen as a block-diagonal matrix; its real eigenvalues are $0$, $0.85$ and $-0.85$ and there are also five pairs of complex eigenvalues, whose arguments are equispaced in a narrow range of frequency where the sinusoids lie.
Firstly, the spectral lines were fixed in $\omega_1=0.42$ rad/s and $\omega_2=0.53$ rad/s, and so the complex poles of $G(z)$ were chosen as:
$$
0.9\rm{e}^{\pm j0.42},\,0.9\rm{e}^{\pm j0.44},\,0.9\rm{e}^{\pm j0.46},\,0.9\rm{e}^{\pm j0.48},\,0.9\rm{e}^{\pm j0.50}.
$$
By considering the constant prior, equal to the sample covariance of the available data, the proposed method was able to detect both lines.
We considered then the more challenging task when  $\omega_1=0.45$ rad/s and $\omega_2=0.47$ rad/s. This choice makes the value of the distance between the two lines lower than the resolution limit of the periodogram, which amounts to $\frac{2\pi}{N}$ (which in our case is  $\frac{2\pi}{300}\simeq 0.021$ rad/s).
Nevertheless, choosing the poles closer to the unit circle, by fixing their radius to $0.95$, the RER estimator was still able to detect the presence of two lines. Figure \ref{fig:testRigheVicine} compares its performances with those achieved by THREE. In simulations, RER exhibited performances that were quite similar or slightly better than those of THREE, as in the case which is shown, where the peaks that were estimated are slightly closer to the real position of the spectral lines than those indicated by THREE. It was also observed that, in general, bringing the poles closer to the unit circle increases both the resolution and the variance of the estimates. The same trade-off was first described in \cite{BYRNES_GEORGIOU_LINDQUIST_THREE} and seems to be typical of all THREE-like methods.

\subsection{Multivariate Case}
In order to test the performances of the proposed method in multivariate spectral estimation, we considered the same estimation task that is described in \cite[Section VIII.C]{RAMPONI_FERRANTE_PAVON_GLOBALLYCONVERGENT}.
The process $y$ was obtained by filtering a bivariate Gaussian white noise process with zero mean and variance $I$ through a square shaping filter of order $40$. The filter coefficients were chosen at random, except for one fixed complex poles pair, $0.9\rm{e}^{\pm j0.52}$ and the zeros pair $(1-10^{-5})\rm{e}^{\pm j0.2}$.

\begin{figure*}[h]
\centering
\includegraphics[width=0.78\textwidth]{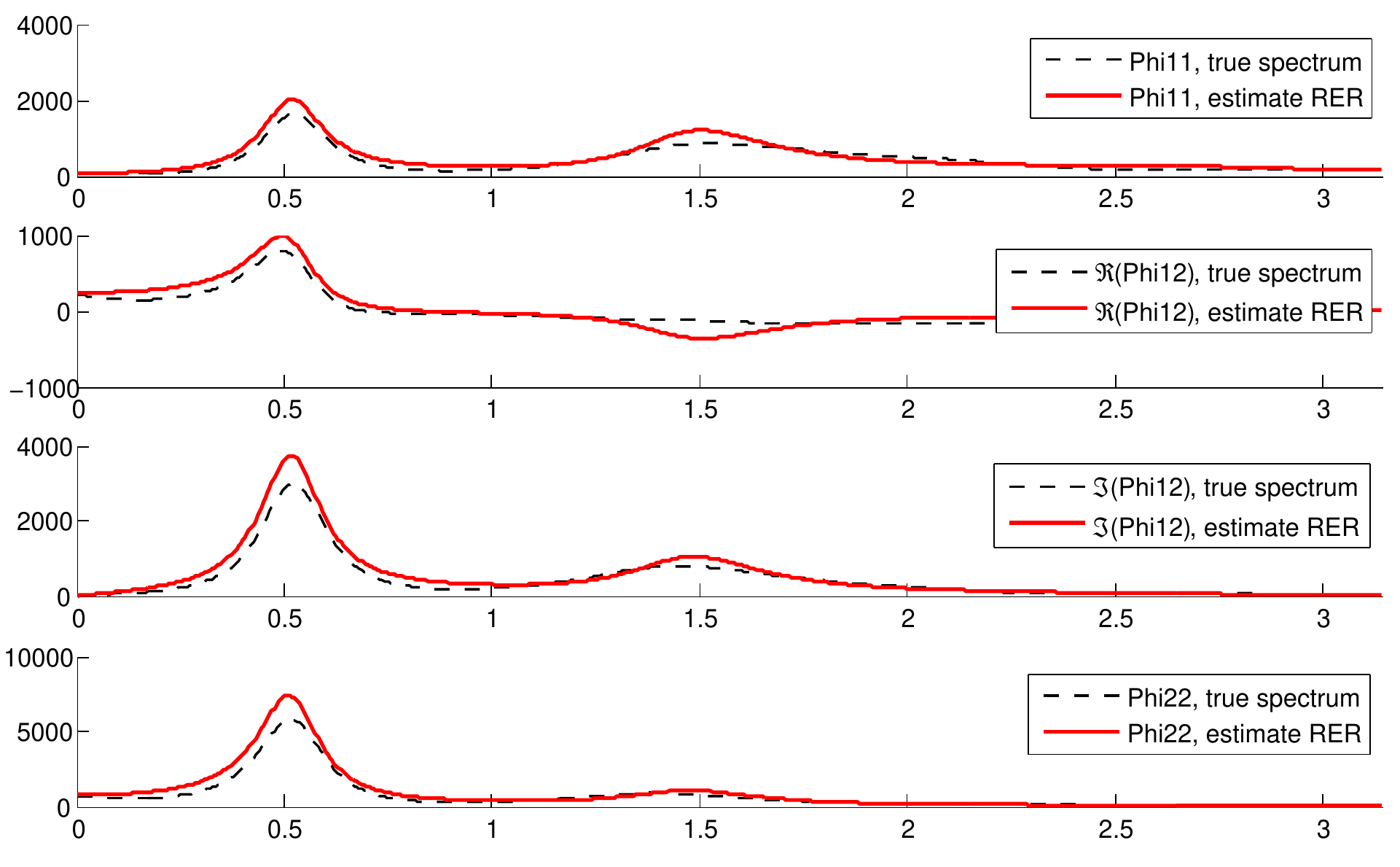}
\caption{Multivariate spectrum estimation, $N=300$, PEM(3) prior. Comparison between the approximant and the true spectrum.}
\label{fig:mvSpectra}
\end{figure*}

\begin{figure}[h]
\centering
\includegraphics[width=0.80\textwidth]{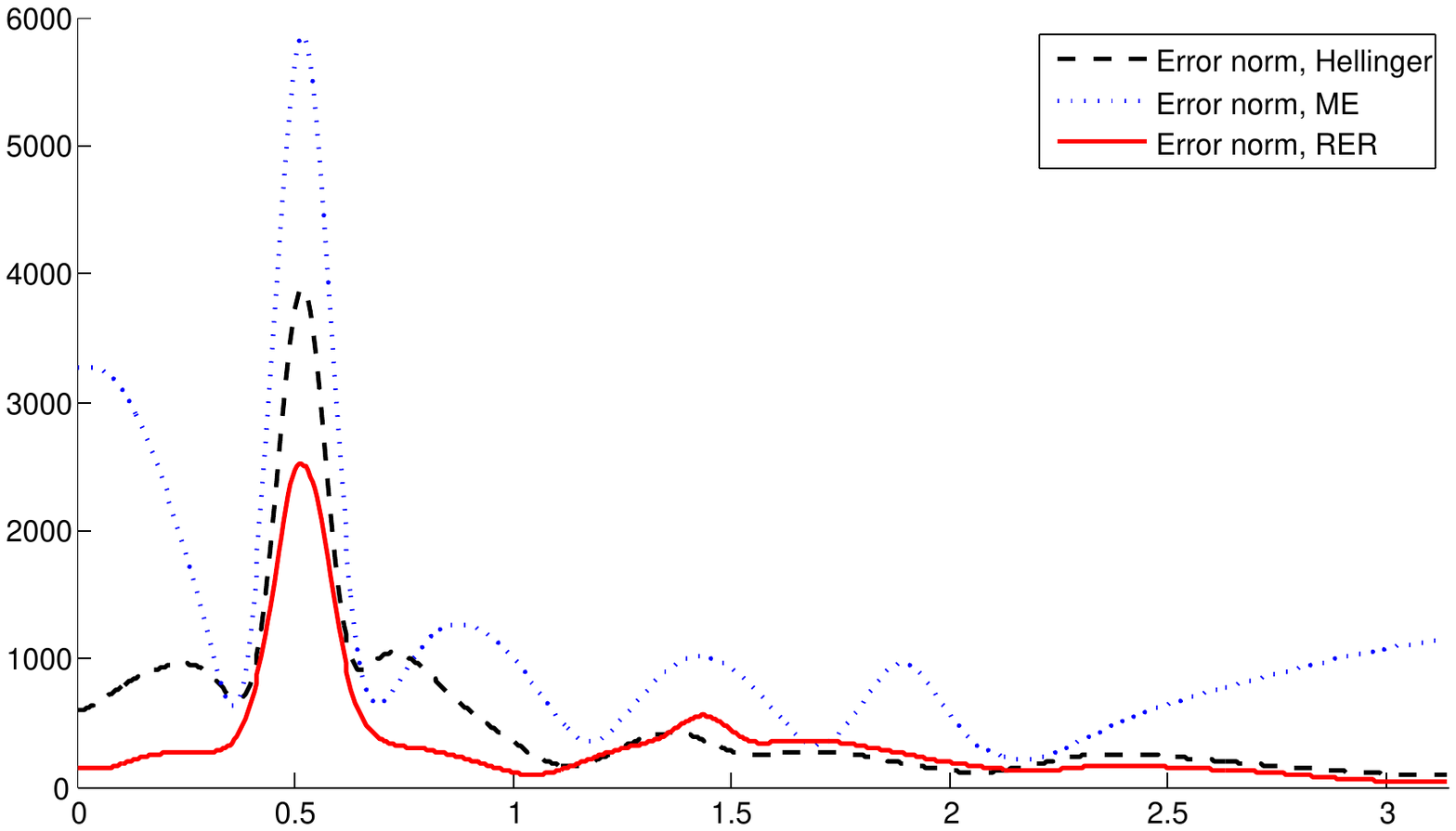}
\caption{ Comparison of THREE-like approaches in terms of average estimation error. $N=300$ available data. Both RER and Hellinger estimator are provided with a PEM(3) prior. All the considered methods make use of the same filter $G(z)$}
\label{fig:mvErrorNorm}
\end{figure}

We designed the filter $G(z)$ by fixing four complex poles pairs with radius $0.7$ and arguments equispaced in the range $\left[0,\pi\right]$. We assumed $N=300$ samples of the process $\{y_k\}_{k\in\mathbb{Z}}$ to be available. As for the the prior, our choice was to compute a simple PEM model of order $3$, by means of the standard function \verb* pem  provided in \textsc{Matlab}.
Figure \ref{fig:mvSpectra} shows the real spectrum and the estimate computed by the RER approach.
We then compared the performance of the proposed technique to those achieved by Maximum Entropy \cite{GEORGIOU_MAXIMUMENTROPY} and Hellinger-distance estimators, which are both THREE-like approaches to multivariate spectral estimation.
Notice that the latter represented the state-of-the art of the convex optimization approach to spectral estimation in the multichannel framework \cite{FERRANTE_PAVON_RAMPONI_HELLINGERVSKULLBACK}.
In order to make the comparison as independent as possible of the specific data set, we performed $50$ trials by feeding the shaping filter with independent realizations of the input noise process.
The criterion to evaluate the performances of each method was taken to be the average estimation error at each frequency, defined as
\eq
E_M(\vartheta):= \frac{1}{50}\sum_{i=1}^{50}\|\hat{\Phi}_{M}(\ejth)-\Phi(\ejth)\|.
\eeq
Here M denotes the specific algorithm, $\hat{\Phi}_{M}$ is the corresponding approximant and the norm is the {\em spectral norm} (i.e. the largest singular value).
Figure \ref{fig:mvErrorNorm} allows to compare the various techniques.
The results achieved by our approach are quite better than those of Maximum Entropy estimator (referred to as ME). Our RER method seems also to slightly outperform the Hellinger-distance approach. It is worth noticing that in the Hellinger case the order of the estimates was $19$, while in the case of RER it was just $11$. The order of the ME estimate  is equal to $8$.

\begin{figure}
\centering
\includegraphics[width=0.80\textwidth]{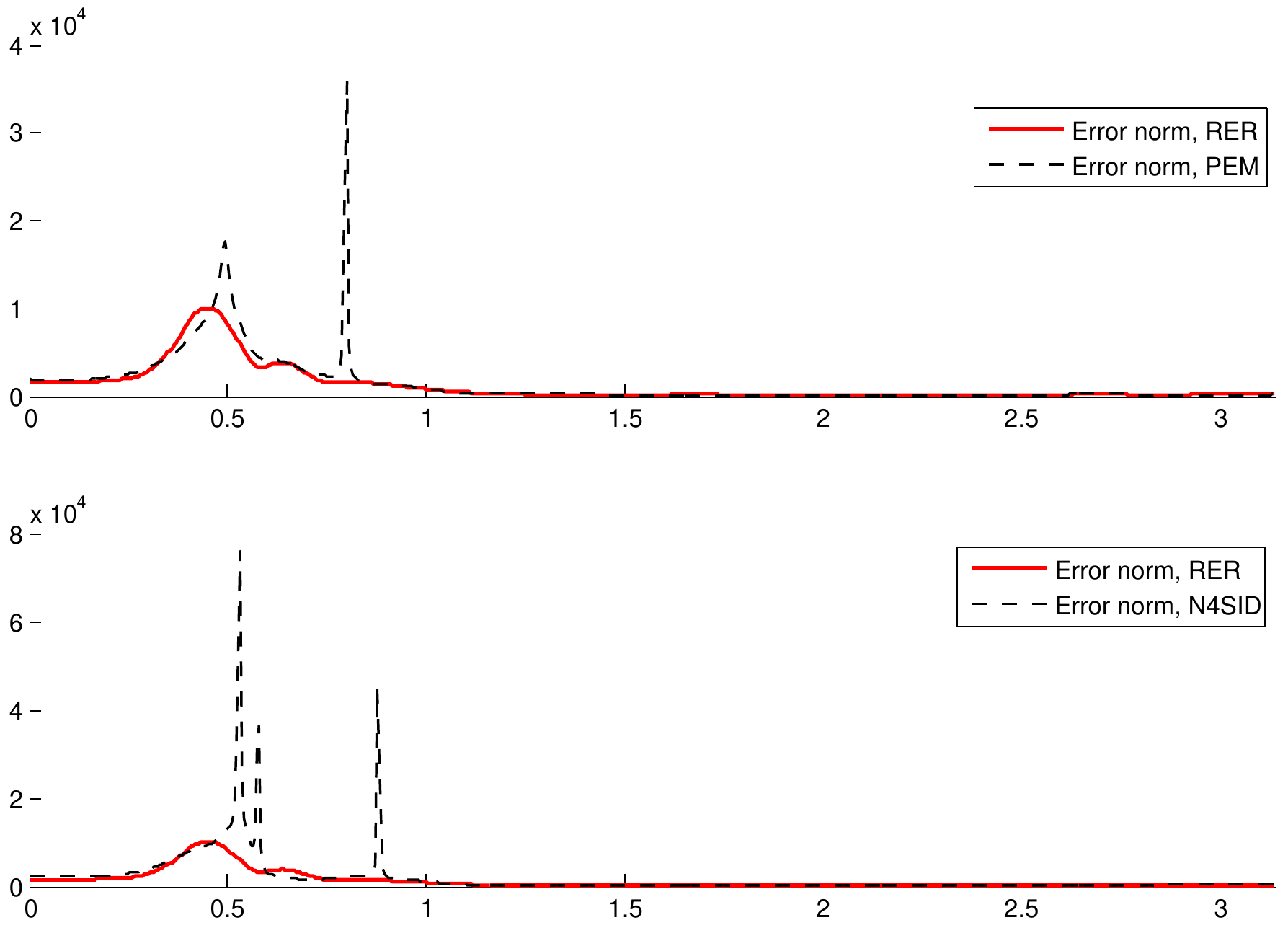}
\caption{Comparison of \textsc{Matlab}'s PEM, \textsc{Matlab}'s N4SID and RER in terms of average estimation error. $N=100$ available samples. The prior considered for RER is a PEM(2) model. The filter $G(z)$ has a pole in the origin and four complex conjugate poles pairs with radius $0.7$. Notice that RER does not exhibit  artifacts, whereas PEM and N4SID do.}
\label{fig:artifacts}
\end{figure}

{ It is interesting to investigate the case when only a few samples of the process $y$ are available.
Shortness of the available data record, can heavily affect with {\em artifacts} the estimates obtained by classical methods such as \textsc{Matlab}'s PEM and N4SID.
As the other THREE-like approaches, the RER method seems to be quite robust with respect to such a problem.     Figure \ref{fig:artifacts} shows the results   obtained in a case where only $N=100$ samples are available.
Both PEM and N4SID estimates are affected by artifacts. On the contrary, the proposed approach is not. This result seems to suggest that RER estimation is suitable to tackle spectral estimation issues characterized} by the presence of short data records of the process of interest.
\FloatBarrier
\section{On the spectral representation of stationary Gaussian processes}\label{spectralprelim}

{We { shall} need the following result whose proof is given in { the} Appendix.
\begin{lemm}\label{strtdipvia}
Let $u,v$ be $k$-dimensional, real random vectors with probability densities $p,q$, respectively.
Let $f:\R^k\rightarrow \R^h$ be measurable and $p_a,q_a$ be the probability densities of the { augmented} vectors $[u^\top\ f(u)^\top]^\top$ and $[v^\top\ f(v)^\top]^\top$, respectively.
Then
\begin{equation}\D(p_a\|q_a)=\D(p\|q).
\end{equation}
\end{lemm}}
\noindent
We { also} need to consider zero-mean, $n$-dimensional, complex-valued Gaussian random vectors $z=\alpha+j\beta$, where the real and imaginary parts are jointly Gaussian.
The corresponding density function is  the joint probability density of the $2n$-dimensional compound vector $\gamma=[\alpha\tp \beta\tp]\tp$.
The differential entropy of the $n$-dimensional complex Gaussian { density} $p$ { (defined exactly as in (\ref{DiffEntropy}) but with integration taking place over $\R^{2n}$)}  with zero mean is given by
\begin{equation}\label{DiffEntropycomplex}
H(p)=-\int_{\R^{2n}}\log (p(x))p(x)dx=\frac{1}{2}\log(\det R)+\frac{1}{2}(2n)\left(1+\log(2\pi)\right),
\end{equation}
where $R$ is the covariance matrix of the $2n$-dimensional vector $\gamma$.
Similarly, the relative entropy between two zero-mean $n$-dimensional complex Gaussian { densities} $p$ and $q$ is given by
\begin{equation}\label{divgausscomplex}
\D(p\|q):=\frac{1}{2}\left[ \log\det
(R_p^{-1}R_q)+\tr(R_q^{-1}R_p)-2n \right],
\end{equation}
where $R_p$ and $R_q$ are the covariance matrices of the $2n$-dimensional vectors $\gamma_p$ and { $\gamma_q$} corresponding to the { densities} $p$ and $q$, respectively.
If the zero-mean, $\C^n$-valued Gaussian random vector $z$ has the property that
$\E[zz^\top]=0$,   then we say that $z$ is a {\em circular symmetric}
 normally distributed random vector \cite{MillerCSP}. { This implies} that $\E[\alpha\alpha^T]=\E[\beta\beta^T]$.
If $p$ and $q$ are two $n$-dimensional complex Gaussian distribution with circular symmetry, the expression of the relative entropy simplifies to the formula
\begin{equation}\label{divgausscomp_circ_symm}
\D(p\|q)= \log\det (P^{-1}Q)+\tr(Q^{-1}P)-n,
\end{equation}
where  $P$ and $Q$ are the covariance matrices of $\gamma_p$ and $\gamma_q$, respectively.

We now state a few basic facts about the spectral representation of a stationary process that can be found, for instance, in \cite{Kramer-Leadbetter,ROZANOV_STATIONARY,LPBook}.
Let $y=\{y_k;\,k\in\Z\}$ be a $\R^m$-valued, zero-mean, Gaussian, stationary  process and let $C_l:=E\{y_{k+l}y_k\trasp \},\, l\in\Z$, be its {\em covariance lags}. Then
\begin{equation}\label{covspectral}
C_l=\frac{1}{2\pi}{ \int_{-\pi}^\pi}e^{{\imunit} l\vartheta}{\dF},
\end{equation}
where $F$ is a bounded, non-negative, matrix-valued measure   called {\em spectral measure}. The stationary process $y$ admits itself the spectral representation
\begin{equation}
y_k={ \int_{-\pi}^\pi}e^{{\imunit} k\vartheta}d\hat{y}(e^{{\imunit}\vartheta}),
\end{equation}
where $\hat{y}$ is a { $m$-dimensional { stochastic orthogonal measure, see \cite{ROZANOV_STATIONARY}}}.
It may be obtained by defining, as in \cite[pag. 44]{LPBook},
{
\begin{equation}
\chi_k(\vartheta_1,\vartheta_2):=
\begin{cases}
\frac{e^{-{\imunit}\vartheta_2k}-e^{-{\imunit}\vartheta_1k}}{-2\pi {\imunit} k} & \text{if $k\neq 0$}\\
\frac{\vartheta_2-\vartheta_1}{2\pi}& \text{if $k= 0$}
\end{cases}\,,
\end{equation}
}
and setting
\begin{equation}\label{defyhat}
\hat{y}(e^{{\imunit}\vartheta_1},e^{{\imunit}\vartheta_2}):=\lim_{N\rightarrow +\infty}
\sum_{k=-N}^N\chi_k(\vartheta_1,\vartheta_2)y_k
\end{equation}
where the sequence converges in mean square.
We use the notation $d\hat{y}(e^{{\imunit}\vartheta})$ as a short-hand for $\hat{y}(e^{{\imunit}\vartheta},e^{{\imunit}(\vartheta+d\vartheta)})$ (with $d\vartheta>0$). It is well known that\begin{equation}\label{COV}\E\left\{d\hat{y}(e^{{\imunit}\vartheta})d\hat{y}(e^{{\imunit}\vartheta})^*\right\}=\dF,
\end{equation}
where $*$ denotes transposition plus conjugation. If the process $y$ is purely nondeterministic, then $\dF=\Phi_y(e^{{\imunit}\vartheta})d\vartheta$, where $\Phi_y$ is the spectral density function.
\begin{propo}\label{circsymmnormal-indep}
 Suppose $\vartheta_1,\vartheta_2\in(-\pi,\pi]$, then
 $\hat{y}(e^{-{\imunit}\vartheta_2},e^{-{\imunit}\vartheta_1})=\overline{\hat{y}(e^{{\imunit}\vartheta_1},e^{{\imunit}\vartheta_2})}
$. If, moreover, $\vartheta_1,\vartheta_2$ have the same sign, then, $\hat{y}(e^{{\imunit}\vartheta_1},e^{{\imunit}\vartheta_2})$ is  a circularly symmetric, normally distributed random vector.
Finally, let $\vartheta_1,\vartheta_2,\vartheta_3,\vartheta_4$ be such that
 $[\vartheta_1,\vartheta_2] \cap [\vartheta_3,\vartheta_4]=[\vartheta_1,\vartheta_2] \cap [-\vartheta_4,-\vartheta_3]=\emptyset$. Then,  $\hat{y}(e^{{\imunit}\vartheta_1},e^{{\imunit}\vartheta_2})$ and $\hat{y}(e^{{\imunit}\vartheta_3},e^{{\imunit}\vartheta_4})$ are independent random vectors.
 \end{propo}
\begin{IEEEproof}
{Observe that $\hat{y}(e^{{\imunit}\vartheta_1},e^{{\imunit}\vartheta_2})$ is a complex-valued random vector that may be written as $\hat{y}(e^{{\imunit}\vartheta_1},e^{{\imunit}\vartheta_2})=\hat{y}_r(e^{{\imunit}\vartheta_1},e^{{\imunit}\vartheta_2})+\imunit\hat{y}_i(e^{{\imunit}\vartheta_1},e^{{\imunit}\vartheta_2})$. In view of (\ref{defyhat}) the real
part $\hat{y}_r(e^{{\imunit}\vartheta_1},e^{{\imunit}\vartheta_2})$  and the imaginary part
$\hat{y}_i(e^{{\imunit}\vartheta_1},e^{{\imunit}\vartheta_2})$  are jointly Gaussian real random vectors, so $\hat{y}(e^{{\imunit}\vartheta_1},e^{{\imunit}\vartheta_2})$ is a complex-valued Gaussian vector.
Since $y_k$ is a $\R^m$-valued process, $\hat{y}(e^{{\imunit}\vartheta_1},e^{{\imunit}\vartheta_2})$ (which may be thought of as an integrated version of a ``Fourier transform'') has the Hermitian symmetry or equivalently $\hat{y}(e^{-{\imunit}\vartheta_2},e^{-{\imunit}\vartheta_1})=\overline{\hat{y}(e^{{\imunit}\vartheta_1},e^{{\imunit}\vartheta_2})}
$.
Moreover, for $\vartheta_1$ and $\vartheta_2$ with the same sign, $\hat{y}(e^{{\imunit}\vartheta_1},e^{{\imunit}\vartheta_2})$ and $\hat{y}(e^{-{\imunit}\vartheta_2},e^{-{\imunit}\vartheta_1})$ are orthogonal. Thus, we get
$$
0=\E[\hat{y}(e^{{\imunit}\vartheta_1},e^{{\imunit}\vartheta_2})\hat{y}(e^{-{\imunit}\vartheta_2},e^{-{\imunit}\vartheta_1})^*]=\E[\hat{y}(e^{{\imunit}\vartheta_1},e^{{\imunit}\vartheta_2})\hat{y}(e^{{\imunit}\vartheta_1},e^{{\imunit}\vartheta_2})^\top]
$$
or, equivalently, $\hat{y}(e^{{\imunit}\vartheta_1},e^{{\imunit}\vartheta_2})$ is circularly symmetric normally distributed.
Finally, recall that two {\em complex} Gaussian random vectors $v_1$, $v_2$ are independent if and only if  $E[v_1v_2^\top]=E[v_1v_2^\ast]=0$.
In our case, by the orthogonality property, we have $E[\hat{y}(e^{{\imunit}\vartheta_1},e^{{\imunit}\vartheta_2})\hat{y}(e^{{\imunit}\vartheta_3},e^{{\imunit}\vartheta_4})^\ast]=0$ and
$E[\hat{y}(e^{{\imunit}\vartheta_1},e^{{\imunit}\vartheta_2})\hat{y}(e^{{\imunit}\vartheta_3},e^{{\imunit}\vartheta_4})^\top]=E[\hat{y}(e^{{\imunit}\vartheta_1},e^{{\imunit}\vartheta_2})\hat{y}(e^{-{\imunit}\vartheta_4},e^{-{\imunit}\vartheta_3})^\ast]=0$.

}
\end{IEEEproof}

\section{Spectral relative entropy rate}\label{SRER}
Consider two zero-mean, jointly Gaussian, stationary, purely nondeterministic stochastic processes $y=\{y_k;\,k\in\Z\}$ and $z=\{z_k;\,k\in\Z\}$ taking values in $\R^m$ with spectral representation
\begin{align}\label{SR1}
y_k={ \int_{-\pi}^\pi}e^{\imunit k\vartheta}d\hat{y}(e^{{\imunit}\vartheta}),&\qquad\E\left\{d\hat{y}(e^{{\imunit}\vartheta})d\hat{y}(e^{{\imunit}\vartheta})^*\right\}=\Phi_y(e^{{\imunit}\vartheta})d\vartheta,\\
z_k={ \int_{-\pi}^\pi}e^{\imunit k\vartheta}d\hat{z}(e^{{\imunit}\vartheta}),&\qquad\E\left\{d\hat{z}(e^{{\imunit}\vartheta})d\hat{z}(e^{{\imunit}\vartheta})^*\right\}=\Phi_z(e^{{\imunit}\vartheta})d\vartheta.\label{SR2}
\end{align}
 Let $\vartheta_k=\frac{\pi k}{n}$,  and  { consider the complex Gaussian random vectors $\deltayk:=\hat{y}(e^{{\imunit}\vartheta_k},e^{{\imunit}\vartheta_{k+1}})$ and {$\deltazk:=\hat{z}(e^{{\imunit}\vartheta_k},e^{{\imunit}\vartheta_{k+1}})$}, with $k=0,1\dots,2n$. Define now the random vectors
\begin{align}
\hat{Y}_k := \bbm \deltayka{0}\\ \vdots\\  \deltayka{k-1}\ebm, \quad
\hat{Z}_k := \bbm \deltazka{0}\\ \vdots\\  \deltazka{k-1}\ebm,\quad k=1,\dots,2n,
\end{align}
and denote their joint probability density by $p(\hat{Y}_{k})$ and $p(\hat{Z}_{k})$, respectively.
\begin{definition}
The {\em spectral relative entropy rate} between between $y$ and $z$ is defined by
the following limit, provided it exists:
\begin{equation}\label{defser}
 {{\dr{d\hat{y}}{d\hat{z}}}:=}\lim_{n\rightarrow\infty}\frac{1}{2n}\D\left(p(\hat{Y}_{2n})\|p(\hat{Z}_{2n})\right).
\end{equation}
 \end{definition}}
We now establish a remarkable connection between time-domain and spectral-domain relative entropy rates.
\begin{theorem} \label{equalrer}Let $y$ and $z$ be as above. Assume that both $\Phi_y$ and $\Phi_z$ are piecewise continuous, coercive spectral densities. The following { equality} holds:
\begin{equation}
\dr{y}{z}= \dr{d\hat{y}}{d\hat{z}}.\label{tworelentropy}
\end{equation}
\end{theorem}

\begin{IEEEproof}
In view of Proposition \ref{circsymmnormal-indep}, the last $n$ components of  $\hat{Y}_{2n}$ are functions (the complex conjugate) of the first $n$ and the same holds for $\hat{Z}_{2n}$.
Hence, in view of Lemma \ref{strtdipvia}, we have $\D\left(p(\hat{Y}_{2n})\|p(\hat{Z}_{2n})\right)=\D\left(p(\hat{Y}_{n})\|p(\hat{Z}_{n})\right)$.
Using again Proposition \ref{circsymmnormal-indep}, we have that the elements of $\hat{Y}_{n}$ are independent random vectors and the same holds for the elements of $\hat{Z}_{n}$. Hence, we have the following additive decomposition:
\be\label{relentropydecomp}
\D\left(p(\hat{Y}_{2n})\|p(\hat{Z}_{2n})\right)=\D\left(p(\hat{Y}_{n})\|p(\hat{Z}_{n})\right) = \sum_{k=0}^{n-1} \DDdue{p(\delta\hat{y}_k)}{p(\delta\hat{z}_k)},
\ee
{
  with $p(\deltayk)$ and $p(\deltazk)$ being the probability densities of the random vector ${ \deltayk}=\hat{y}(e^{\imunit\vartheta_k},e^{\imunit\vartheta_{k+1}})$ and ${\deltazk=}\hat{z}(e^{\imunit\vartheta_k},e^{\imunit\vartheta_{k+1}})$, respectively.
Since $\deltayk$ and $\deltazk$ are jointly Gaussian and circularly symmetric, by (\ref{divgausscomp_circ_symm}) and (\ref{SR1})-(\ref{SR2}), we get,
\begin{equation}
\D\left(p(\deltayk) \|p(\deltazk)\right)=
 \log\det\left[Q_y^{-1}(\vartheta_k,\vartheta_{k+1})Q_z(\vartheta_k,\vartheta_{k+1})\right] \\
 +\tr\left[Q_z^{-1}(\vartheta_k,\vartheta_{k+1})Q_y(\vartheta_k,\vartheta_{k+1})\right]-m,
\label{relentr}
\end{equation}
where, by virtue of the orthogonal increments property,
\eq
Q_y(\vartheta_k,\vartheta_{k+1}):=\int_{\vartheta_k}^{\vartheta_{k+1}}\Phi_y(e^{{\imunit}\xi})d\xi,
\qquad
Q_z(\vartheta_k,\vartheta_{k+1}):=\int_{\vartheta_k}^{\vartheta_{k+1}}\Phi_z(e^{{\imunit}\xi})d\xi.
\eeq
By piecewise continuity and the mean value theorem, we have that, except for a finite number of $k$'s,
 \begin{equation}
\begin{split}
 \D\left(p(\delta\hat{y}_{k}) \|p(\delta\hat{z}_{k})\right)&=\log\det\left[\left(\Phi_y(e^{{\imunit}{\bar{\vartheta}_k}})\frac{\pi}{n}\right)^{-1}\Phi_z(e^{{\imunit}{\bar{\vartheta}_k}})\frac{\pi}{n}\right]+\tr\left[\left(\Phi_z(e^{{\imunit}{\bar{\vartheta}_k}})\frac{\pi}{n}\right)^{-1}\Phi_y(e^{{\imunit}{\bar{\vartheta}_k}})\frac{\pi}{n}\right]-m\\
 &=\log\det[\Phi_y(e^{{\imunit}{\bar{\vartheta}_k}})^{-1}\Phi_z(e^{{\imunit}{\bar{\vartheta}_k}})]+\tr\left[\Phi_z(e^{{\imunit}{\bar{\vartheta}_k}})^{-1}\Phi_y(e^{{\imunit}{\bar{\vartheta}_k}})\right]-m,
\end{split}
\end{equation}
where $\vartheta_k\leq \bar{\vartheta}_k <  \vartheta_{k+1}$. By employing the latter expression together with \eqref{relentropydecomp} and (\ref{defser}), we get
\begin{equation*}
\begin{split}
\dr{d\hat{y}}{d\hat{z}}&=\lim_{n\rightarrow\infty}\frac{1}{2n}\DDdue{\hat{Y}_n}{\hat{Z}_n}=\lim_{n\rightarrow\infty}\frac{1}{2n}\sum_{k=0}^{n-1} \D\left(p(\delta\hat{y}_{k}) \|p(\delta\hat{z}_{k})\right) \\
&=\lim_{n\rightarrow\infty}\frac{1}{2n}\sum_{k=0}^{n-1} \log\det\Phi_y(e^{{\imunit}{\bar{\vartheta}_k}})^{-1}\Phi_z(e^{{\imunit}{\bar{\vartheta}_k}})
+\tr\left[\Phi_z^{-1}(e^{{\imunit}{\bar{\vartheta}_k}})\left(\Phi_y(e^{{\imunit}{\bar{\vartheta}_k}})-\Phi_z(e^{{\imunit}{\bar{\vartheta}_k}})\right)\right] \\
&=\lim_{n\rightarrow\infty}\frac{1}{2\pi}\sum_{k=0}^{n-1} \left\{\log\det\Phi_y(e^{{\imunit}{\bar{\vartheta}_k}})^{-1}\Phi_z(e^{{\imunit}{\bar{\vartheta}_k}})+\tr\left[\Phi_z^{-1}(e^{{\imunit}{\bar{\vartheta}_k}})\left(\Phi_y(e^{{\imunit}{\bar{\vartheta}_k}})-\Phi_z(e^{{\imunit}{\bar{\vartheta}_k}})\right)\right]\right\}\frac{\pi}{n}\\
&= \frac{1}{2\pi}\int_{0}^\pi \left\{\log\det\left(\Phi_y^{-1}(e^{{\imunit}\vartheta})\Phi_z(e^{{\imunit}\vartheta})\right)+\tr\left[\Phi_z^{-1}(e^{{\imunit}\vartheta})\left(\Phi_y(e^{{\imunit}\vartheta})-\Phi_z(e^{{\imunit}\vartheta})\right)\right]\right\}d\vartheta\\
&= \frac{1}{4\pi}{\int_{-\pi}^\pi} \left\{\log\det\left(\Phi_y^{-1}(e^{{\imunit}\vartheta})\Phi_z(e^{{\imunit}\vartheta})\right)+\tr\left[\Phi_z^{-1}(e^{{\imunit}\vartheta})\left(\Phi_y(e^{{\imunit}\vartheta})-\Phi_z(e^{{\imunit}\vartheta})\right)\right]\right\}d\vartheta,
\end{split}
\end{equation*}
which, by (\ref{relentropyconnection}), is (\ref{tworelentropy}).}
\end{IEEEproof}

\begin{remark}
As is well known, the fundamental property of the Fourier transform is that it is isometric. The above result may be interpreted as a  further  invariance principle of the Fourier transform: the relative entropy rate is the same in the time and spectral domain.
\end{remark}

\section{Conclusion}
{ In this paper, a profound information-theoretic result relating time and spectral domain relative entropy rates of stationary Gaussian processes has been established.
Motivated by this result, a new THREE-like approach to multivariate spectral
estimation, called RER, has been introduced and tested. It  appears as
the most natural extension of maximum entropy methods when a prior
estimate of the spectrum is available. It features an upper bound on
the complexity of the estimate which is equal to the one provided by
THREE in the scalar context, sensibly improving on the best one so far
available in the multichannel setting with prior estimate.} As for
previous THREE-like methods, RER exhibits high resolution features
and works extremely well with short observation records
outperforming \textsc{Matlab}'s PEM and \textsc{Matlab}'s N4SID.
\section*{Acknowledgments}
We wish to thank prof. Paolo Dai Pra for providing the proof of Lemma \ref{strtdipvia}. The constructive comments of four anonymous reviewers are also gratefully acknowledged.

\appendix
\begin{IEEEproof}[Proof of Lemma \ref{strtdipvia} \cite{DAIPRA}]
Recall the variational formula for relative entropy \cite{Dembo1993}:
\begin{equation}\label{eq:var_relentropy}
\D(p\|q)=\sup_{\varphi\,\in\,\phi}\left\lbrace{\E{\left[\varphi(v)\right]} -
\log\E{\left[\e^{\varphi(u)}\right]}}\right\rbrace,
\end{equation}
where $\phi$ is the set of all measurable and bounded functions
$\varphi\,:\,\mathbb{R}^k\rightarrow\mathbb{R}$.
Consider a measurable and bounded function $\varphi\,:\,\R^k\rightarrow\R$. Define $\varphi_a\,:\,\R^{k+h}\rightarrow\R$ by
\be
\varphi_a([x^\top\ x'^\top]^\top):=\varphi(x),
\ee
where $x'\,\in\,\R^h$.
Obviously, $\varphi_a$ is bounded and measurable, and
\be
\E{\left[\varphi(v)\right]}-\log\E{\left[\e^{\varphi(u)}\right]}=\E{\left[\varphi_a(v_a)\right]}-\log\E{\left[\e^{\varphi_a(u_a)}\right]}\leq \D{(p_a\|q_a)}.
\ee
By { taking} the supremum, { we get} that $\D{(p\|q)}\leq\D{(p_a\|q_a)}$.
The opposite inequality can be proven along the same lines. Indeed, let $\psi_a\,:\, \R^{h+k}\rightarrow \R$ be a measurable and bounded function. Define $\psi\,:\,\R^k\rightarrow\R$ by $\psi(x):=\psi_{ a}(x,f(x))$. Then, $\psi$ is measurable and bounded too, so that
\be
\E{\left[\psi_a(v_a)\right]}-\log\E{\left[\e^{\psi_a(u_a)}\right]}=\E{\left[\psi(v)\right]}-\log\E{\left[\e^{\psi(u)}\right]}\leq\D{(p\|q)}.
\ee
In view of \eqref{eq:var_relentropy}, we now get $\D{(p_a\|q_a)}\leq\D{(p\|q)}$.
\end{IEEEproof}

\begin{IEEEproof}[Proof of Lemma \ref{lemma:bd}]
\begin{enumerate}
\item
As a consequence of Lemma \ref{lemm:tracciaLambda},
\begin{equation}
 \begin{split}
  J_\Psi^n(\Lambda)&= \int \tr \left[ \Lambda - \log(I+G_1^* \Lambda G_1 + \frac{1}{n}I)\right]\\
                 &\geq  \int \tr \left[\mu (I+G_1^* \Lambda G_1) - \log(I+G_1^* \Lambda G_1 + \frac{1}{n}I)\right]+\alpha.
   \end{split}
\end{equation}
Let $\left\lbrace x_i \right\rbrace$ be the eigenvalues of $(I+G_1^* \Lambda G_1)$. Then,
\begin{equation}
 \begin{split}
  J_{\Psi}^n(\Lambda)&=\int \tr \left[\mu (I+G_1^* \Lambda G_1) - \log(I+G_1^* \Lambda G_1 + \frac{1}{n}I)\right]+\alpha\\
                       &=\int \mu\sum_{i=1}^m x_i - \sum_{i=1}^m \log\left(x_i + \frac{1}{n}\right)+\alpha=\int \rho\left(x_1,\dots,x_m\right)+\alpha,
 \end{split}
\end{equation}
where $\rho(x_1,\dots,x_m):=\mu\sum_{i=1}^m x_i - \sum_{i=1}^m \log\left(x_i + \frac{1}{n}\right)$.
Moreover,
$$
\frac{\partial }{\partial x_i} \left[\rho(x_1,\dots,x_m)\right] = \mu - \frac{1}{x_i + \frac{1}{n}}\quad \forall \,i.
$$
The minimum of $\rho$ is thus attained by choosing $x_i=\frac{1}{\mu}- \frac{1}{n}$, $\forall \, i$. Therefore,
$$
\rho(x_1,\dots,x_m)\geq m -\frac{\mu m}{n}+m\log{\mu}
$$
The fact  that $J_\Psi^n(\Lambda)$ is { bounded from below} over $\CL$ now follows:
\begin{equation}
  J_\Psi^n(\Lambda) \geq \alpha +  2\pi[m +m\log{\mu} - \frac{\mu m}{n}]  \geq \alpha + 2\pi m[1 +\log{\mu}].
\end{equation}
\item Beppo Levi's Theorem allows to conclude that $J_\Psi^\infty(\Lambda)=J_\Psi(\Lambda)$ in $\LL$:
\begin{equation}
 J_\Psi^\infty(\Lambda)= \int \tr[\Lambda] - \int \tr \left[\lim_{n\rightarrow\infty} \log(I+G_1^*\Lambda G_1 + \frac{1}{n}I)\right] = J_\Psi(\Lambda).
\end{equation}
\item Since, for $\Lambda \,\in \,\mathcal{B}^c$, the rational function $\det{(I+G_1^* \Lambda G_1)}$ is not identically zero, its logarithm is integrable over {$(-\pi,\pi]$}. Hence, $J_\Psi^\infty(\Lambda)$ is finite.  $J_\Psi^\infty(\Lambda)=+\infty$ instead for $\Lambda \, \in \, \mathcal{B}$.
\end{enumerate}
\end{IEEEproof}

\begin{IEEEproof}[Proof of Lemma \ref{lemma:infJ}]
In view of Lemma \ref{lemm:tracciaLambda}
\begin{equation}\label{eq:bobel}
\tr\left[\Lambda\right]\geq \mu \tr \left[\int (G_1^*\Lambda G_1+I)\right] +\alpha>\alpha,
\end{equation}
 so that $\tr\left[\Lambda\right]$ is { bounded from below}.
Consider a sequence $\left\lbrace\Lambda_k\right\rbrace_{k\in\mathbb{N}} \, \in \, \LL$, such that $$\lim_{k \rightarrow \infty} \|\Lambda_k\| = +\infty.$$
Let $\Lambda_k^0:=\frac{\Lambda_k}{\|\Lambda_k\|}$. Since $\LL$ is convex and $\Lambda=0$ belongs to $\LL$, $\forall \, {\xi} \, \in \left[0,1\right]$, ${\xi}\Lambda\in\,\LL$. Therefore $\Lambda_k^0 \, \in \LL$ for sufficiently large $k$.
Let $\eta:=\liminf \tr \left[\Lambda_k^0\right]$
In view of (\ref{eq:bobel})
$$
\tr\Lambda_k^0 = \frac{1}{\|\Lambda_k\|}\tr\Lambda_k > {\frac{1}{\|\Lambda_k\|}} \alpha \, \rightarrow \, 0,$$
for $\|\Lambda_k\|\rightarrow \infty$, so $\eta \geq 0$.
Thus, the sequence { $\left\lbrace\Lambda_k^0\right\rbrace$} has
a subsequence such that the limit of its trace is $\eta$. Given that $\Lambda_k^0$ belongs to the surface of the unit ball, which is compact, the subsequence contains a subsubsequence $\left\lbrace\Lambda_{k_m}^0\right\rbrace_{k_m\in\mathbb{N}}$ that is convergent.
Define
$$\Lambda_\infty ^0 := \lim_{k_m \rightarrow \infty} \Lambda_{k_m}^0.$$
The next step is to prove that $\Lambda_\infty ^0\, \in \LL$. To this aim, notice that $\Lambda_\infty^0$ is the limit of a convergent sequence in the finite-dimensional linear space $\Range{(\Gamma)}$. Therefore it belongs to $\Range{(\Gamma)}$. Moreover, recall that the primary sequence $\left\lbrace\Lambda_k\right\rbrace_{k\in\mathbb{N}}$ has elements belonging to $\LL$. It means that, for each $\Lambda_k$, $(I+G_1^*\Lambda_k G_1)>0$. As a consequence, it holds that, for each $m$,
$$\left(\frac{1}{\|\Lambda_{k_m}\|}I+G_1^*\Lambda_{k_m}^0G_1\right)>0 \quad \text{on $\mathbb{T}$ }.$$
Taking the pointwise limit for $m\rightarrow \infty$, it results that $G_1^*\Lambda_\infty^0 G_1$ is positive semidefinite on $\T$, and so $(I+G_1^*\Lambda_\infty^0 G_1)$ is strictly positive definite on $\mathbb{T}$. Therefore, $\Lambda_\infty^0 \, \in \LL$.

The next step is to prove that $\tr \Lambda_\infty^0>0$.
If the feasibility condition (\ref{eq:feasibility}) holds, there exists $\Phi_I$ such that $I=\int G\Phi_I G^*$. Therefore, it is possible to write:
\begin{equation}\label{eq:trPos}
\begin{split}
\tr\Lambda_\infty^0  &= \tr \int G\Phi_I G^*\Lambda_\infty^0 =\int \tr  \left[W_\Psi^{-*}W_\Psi^*G^*\Lambda_\infty^0 G W_\Psi W_\Psi^{-1} \Phi_I\right]\\
                                     &= \int \tr \left[G_1^*\Lambda_\infty^0 G_1 \underbrace{W_\Psi^{-1} \Phi_I W_\Psi^{-*}}_{\Xi}\right]=\int \tr \left[\Xi^{\frac{1}{2}} G_1^*\Lambda_\infty^0 G_1 \Xi^{\frac{1}{2}}\right],
\end{split}
\end{equation}
where the coercive spectral density $\Xi$ is defined as in Lemma \ref{lemm:tracciaLambda}.
Since $G_1^*\Lambda_\infty^0 G_1 \geq 0$, in order to prove that $\tr \left[\Lambda_\infty^0\right]$ is positive, in view of \eqref{eq:trPos} it is sufficient to show that $G_1^*\Lambda_\infty^0 G_1 $ is not identically zero.
Assume by contradiction that $G_1^*\Lambda_\infty^0 G_1 \equiv 0$. As a consequence, $\forall \,\ejth \in \T$,
\begin{equation}
  0 \equiv G_1^*\Lambda_\infty^0 G_1 = W_\Psi^*G^*\Lambda_\infty^0 GW_\Psi.
\end{equation}
Therefore, $G^*\Lambda_\infty^0 G\equiv 0$. However, this means that $\Lambda_\infty^0 \in\, \Range(\Gamma)^\perp$. But it has already been proven that $\Lambda_\infty^0 \in \, \Range(\Gamma)$. Moreover, $\Lambda_\infty^0 \neq 0$, since it belongs to the surface of the unit ball. This is a contradiction. Thus, $G_1^*\Lambda_\infty^0 G_1$ is not identically zero, and from (\ref{eq:trPos}) it follows that $\eta=\tr \Lambda_\infty^0 > 0$.
It follows that there exists $K$ such that $\tr \Lambda_k^0>\frac{\eta}{2}$ for all $k>K$. Notice that $G_1^*G_1$ is positive definite on $\mathbb{T}$ (and indeed coercive). Moreover, $G_1^*\Lambda_k^0 G_1 \leq G_1^*G_1$, since $\Lambda_k^0$ belongs to the unit ball. Therefore,
\begin{equation}\nonumber
 \begin{split}
  \liminf_{k\rightarrow\infty} J_\Psi(\Lambda_k)&= \liminf_{k\rightarrow \infty} \int \tr\left[\Lambda_k-\log(I+G_1^*\Lambda_k G_1)\right]\\
    &=\liminf_{k\rightarrow \infty}  \tr \left[\|\Lambda_k\|\Lambda_k^0\right] -\liminf_{k\rightarrow \infty}\int\tr\left[ \log\left[ \|\Lambda_k\|\left(\frac{1}{\|\Lambda_k\|}I+G_1^*\Lambda_k^0 G_1\right)\right]\right]\\
    &\geq \liminf_{k\rightarrow \infty}  \|\Lambda_k\|\frac{\eta}{2} -\liminf_{k\rightarrow \infty}\int\log\|\Lambda_k\| -\liminf_{k\rightarrow \infty}\int \tr \left[\log\left(\frac{1}{\|\Lambda_k\|}I+G_1^*G_1\right)\right]\\
    &=\liminf_{k\rightarrow \infty} \frac{\eta}{2}\left( \|\Lambda_k\| - \frac{4\pi}{\eta}\log\|\Lambda_k\|\right) -\liminf_{k\rightarrow \infty}\int \tr \left[\log\left(\frac{1}{\|\Lambda_k\|}I+G_1^*G_1\right)\right]\\
    &= +\infty.
\end{split}
\end{equation}
\end{IEEEproof}

\bibliographystyle{plain}
\bibliography{biblio}

\end{document}